\documentclass[12pt]{article}
\usepackage{amsfonts}
\usepackage{amsthm}
\usepackage{amssymb}
\usepackage{amsmath}
\usepackage{graphicx,color,epstopdf}
\usepackage{hyperref}
\def\hDash{\bot\!\!\!\bot}

\newtheorem{theorem}{Theorem}[section]

\newtheorem{lemma}{Lemma}[section]

\numberwithin{equation}{section}

\begin{document}

\title{ Dimensionality determination: a thresholding double ridge ratio criterion
\footnote{Lixing Zhu is a Chair professor of Department of Mathematics
at Hong Kong Baptist University, Hong Kong,  and a professor of School of Statistics at Beijing Normal University, Beijing, China. Xuehu Zhu is an assistant professor at Xi'an Jiaotong University, Xi'an, China and Tao Wang is a research professor at Shanghai Jiaotong University, Shanghai, China. Lixing Zhu's research  was supported by a grant from the
University Grants Council of Hong Kong, Hong Kong, China.}}
\author{Xuehu Zhu, Tao Wang and Lixing Zhu
}
\date{}
\maketitle

\renewcommand\baselinestretch{1.5}
{\small}

\noindent {\bf Abstract.}
Popularly used eigendecomposition-based criteria such as BIC type, ratio estimation  and principal component-based criterion often underdetermine model dimensionality for regressions or the number of factors for factor models. This longstanding problem is caused by the existence of one or two dominating eigenvalues compared to other nonzero eigenvalues. 
 To alleviate this difficulty, we propose a thresholding double ridge ratio criterion such that the true dimension can be better identified and is less underdetermined.  Unlike all existing eigendecomposition-based criteria, this criterion can define consistent estimate without requiring the uniqueness of minimum and can then handle possible multiple local minima scenarios. This generic strategy would be readily applied to other dimensionality or order determination problems. In this paper, w e systematically investigate, for general sufficient dimension reduction theory, the dimensionality determination with fixed and divergent dimensions; for local alternative models that converge to its limiting model with fewer projected covariates, discuss when  the  number of projected covariates can be consistently estimated, when cannot; and for ultra-high dimensional factor models, study  the estimation consistency for the number of common factors.
Numerical studies are conducted to examine the finite sample performance of the method.
\

\bigskip


\newpage
\baselineskip=21pt

\newpage

\setcounter{equation}{0}
\section{Introduction}
The research described herewith is motivated dimensionality determination in three problems in sufficient dimension reduction, model checking for regressions and approximate factor models. As these three problems relate to the eigendecomposition-based criteria for respective target matrices, we begin with the problem in sufficient dimension reduction.

Let $X=(X_1,\ldots,X_p)^{\top}\in \mathbb{R}^{p}$ be   the vector of $p$ covariates,
 $Y=(Y_1,\ldots,Y_m)^{\top}\in \mathbb{R}^{m}$ be the vector of $m$ responses.
If the dependence on $X$ of $Y$ is unspecified, the curses of dimensionality arises and  statistical analysis is difficult even when the number of the covariates  is moderate.
Effective dimension reduction / sufficient dimension reduction (Li, 1991; Cook, 1998) can efficiently reduce the number of the covariates through finding a few most informative linear combinations of the covariates.    The basic idea
is, without loss of information between $Y$ and $X$, to replace the covariates $X$ with its projection onto a $q$-dimensional subspace, assuming $Y$ is linked to $X$ via $q$ linear combinations $B^{\top}X$ of $X$.  Here, $B\in\mathbb{R}^{p \times q}$ and $q$ is usually much smaller than $p$. After such a dimension reduction, variable selection can also be performed to further rule out  unimportant covariates in the linear combinations  such that the models can be more parsimonious (see, e.g. Li et al., 2005; Chen et al., 2010; Wang and Zhu, 2013; Wang et al., 2015).

In this field,   Li (1991) in his seminal paper proposed sliced inverse regression (SIR), see also a relevant reference, Duan and Li (1991).  Cook (1998) subsequently introduced the concept of sufficient dimension reduction (SDR). The success of these methods hinges on successfully identifying the subspace $\emph{S}$ of minimal dimension. Under the central subspace framework,  SDR is to seek a minimum subspace $\emph{S}$ such that
\begin{eqnarray}\label{cen}
Y \hDash X|P_{\emph{S}}X,
\end{eqnarray}
where $\hDash$ denotes the statistical independence and $P_{(\cdot)}$ stands for a projection operator with respect to the standard inner product. The minimum subspace $\emph{S}$
is called the central subspace, and is denoted as  $\emph{S}_{Y|X}$.
As the mean function of a regression model is often more important, then under the central mean subspace framework (Cook and Li, 2002), the  objective of SDR is to find a minimum subspace $\emph{S}$ satisfying
\begin{eqnarray}\label{cen-mean}
Y \hDash E(Y|X)|P_{\emph{S}}X,
\end{eqnarray}
which is called the central mean subspace, and written as  $\emph{S}_{E(Y|X)}$.
The central variance subspace (Zhu and Zhu, 2009) is defined as the minimum subspace $\emph{S}$ that satisfies
\begin{eqnarray}\label{cen-var}
Y \hDash \rm{var}(Y|X)|P_{\emph{S}}X,
\end{eqnarray}
and written as  $\emph{S}_{\rm{var}(Y|X)}$.
Such  minimum dimension reduction subspaces often exist uniquely (Cook, 1998).
Their dimensions are called the structural dimensions of the respective SDR subspaces. In our setup, write the structural dimension as $q$.

Determining the structural dimension plays a vital role in sufficient dimension reduction and most methods are based on eigendecomposition of target matrices.
There are mainly three methodologies available in the literature. The first is hypothesis testing-based methodology  proposed by Li (1991) who developed a sequential testing procedure for sliced inverse regression. Other extensions in this direction include Schott (1994), Velilla (1998), Bura and Cook (2001), 
Cook and Li (2004), Zeng (2008), and Bura and Yang (2011).  
However, hypothesis testing-based methods cannot provide consistent estimates as they  are related to Type-I errors in the sequence of testing procedures, which do not vanish.
Therefore, these methods generally require a relatively large sample size for good performance.
The second is based on information criteria. Zhu et al. (2006) firstly advised a  Bayesian Information Criterion (BIC) type method. There are several modifications in the literature such as Li and Lu (2008) and Wang and Yin (2008). The third methodology uses ratios of eigenvalues to define criteria. Luo et al. (2009) and  Xia et al. (2015) are relevant references.

There are also some other methods that are not eigendecomposition-based. Wang and Xia (2008), following the idea of Xia et al. (2002), proposed a leave-one-out cross-validation (CV) criterion to select $q$ by minimizing the average prediction error. Ma and Zhang (2015),  basing on the estimating equation, suggested a validated information criterion to estimate the structure dimension.
These are unavoidable to  estimate  the nonparametric  regression function with high-dimensional vector of covariates.

We re-visit this issue because of underdetermination that is a longstanding problem in the SDR research field for all aforementioned eigencomposition-based methods. This problem often  causes missing some important linear combinations of covariates and then working models could be too parsimonious to well fit data.  The reason behind is that one or two largest estimated eigenvalues are often dominating and the other estimated eigenvalues are  close to each other in magnitude no matter whether they are nonzero or not at the population level. The minimizer or maximizer of a criterion over all the indexes is often smaller than the true dimension $q$ of the central (mean) subspace.   Zhu et al. (2006) actually had found this phenomenon. See also a relevant reference by Ahn and Horenstein (2013) for estimating the number of common factors in approximate factor models.
 This problem is even more severe  for model checks for regressions (Guo et al., 2015) when we need to identify the central mean subspace of local alternative models. More specifically,
under local alternative hypothesis, we consider a sequence of local alternative models (Guo, et al., 2015):  rewriting a scalar response $Y$ as $Y_n$
\begin{eqnarray}\label{sequence}
Y_n=G(B_1^{\top}X)+C_n g(B^{\top}X)+\varepsilon.
\end{eqnarray}
Here $\varepsilon$ and $X$ are mutually independent, and both the functions $G(\cdot)$ and $g(\cdot)$ are smoothing. $Y_n$ is the response and $X$ is  the vector of covariates independent of $\varepsilon$, $B_1$ is a $p\times q_1$ orthonormal matrix with $q_1\ge 1$, $(B_1, B)$ has a rank $q>q_1\ge 1$ and $C_n$ is a sequence of constants that could depend on the sample size $n$. The model $ Y=G(\beta^{\top}X)+\varepsilon$ is called the hypothetical model (or the limiting model as $C_n\to 0$). When $C_n$ are a fixed constant, it is called the global alternative hypothesis and when $C_n\to 0$, they are called the local alternatives. 
Under certain regularity conditions, the central mean subspace $S_{E(Y_n|X)}$ is equal to  the column subspace spanned by $(B_1, B)$.  Note that the limiting model of this sequence  as $n\to \infty$ is $ Y=G(B_1^{\top}X)+\varepsilon$ such that the corresponding  central mean subspace $S_{E( Y|X)}$ is equal to the space spanned by $B_1$ that is a smaller subspace with the structural dimension $q_1$. However, for any fixed $n$, $C_n\not = 0 $, the central mean subspace $S_{E(Y_n|X)}$  always has a higher structural dimension $q>q_1$. A significant feature of this problem is that
although the local models converge to its limiting model (hypothetical model),  the dimension $q$ does not accordingly shrink to $q_1$, underdetermination is very easy to occur when $C_n$ is too close to $0$. The readers can refer to Theorem~2 of Section~3.2 of Guo et al. (2015) for a special case where $B_1=\beta$ is a vector with $q_1=1$. It is necessary to study at what convergence rate of $C_n\to 0$, $S_{E(Y_n|X)}$ can be identified and at what rate, it cannot.
This problem can be formulated as a special case of general sufficient dimension reduction when the  sequence of the responses depends on the sample size $n$: $Y_n \hDash E(Y|X)|P_{\emph{S}}X$ and its limiting version $ Y \hDash E( Y|X)|P_{\tilde S}X$ where $\tilde S$ is a smaller central mean subspace than $\emph{S}$. The same notion can be extended to the central variance subspace and central subspace. For ease of exposition, we in this paper only focus on the special scenario about model checking for regressions.

This motivates us to develop  a more efficient criterion to determine the dimensionality. The new method has some interesting features as enumerated below.
\begin{enumerate}
\item  `Standardized' eigenvalues are used such that the dominating eigenvalues are not too dominating.
\item At the population level, defined initial ratios of the `standardized' eigenvalues make 
    the $q$-th ratio tends to infinity to be a maximum. 
    This can be a base for a criterion such as the minimizer of the reciprocals of the ratios over all indexes. A relevant reference is Xia et al. (2015). But we find that at the sample level, this is not efficient  as the dominating eigenvalues would still make the first ratio(s) be even smaller than the reciprocal of the $q$th ratio.
\item Second round of ratios is performed to  further distinguish between the  ratio at the true dimension $q$ and the others. The  ratios of the second round  have the following properties. 
    At the population level, the $(q-1)$th ratio tends positive  infinity, the $q$th ratio tends to zero and the others tend to one. The $q$th ratio is then at the valley bottom and the $(q-1)$th is on the peak such that the $q$th ratio can be better separated from the first $q-2$ ratios. This round of ratios can well separate the $q$th ratio from the others. However,   this criterion  would still have more than one local minimum in the previous ratios.  It is particularly the case for local models. Thus, using minimizer to define an estimate would still get a local minimizer that is smaller than $q$.
 \item A largest index, as an estimate, of which the ratio is under a threshold $\tau$ with $0<\tau <1$ can largely avoid multiple local minima problem. The idea of this criterion is unlike all the existing methods that rely on the existence of global minimum/maximum. This method of thresholding would become a new and generic method in this field.
\item To make the ratios of $0/0$ for $j>q$ at population level be stable, ridges in the two round of ratios are added. This is necessary  for constructing ratios.
    \end{enumerate}
We then call the new criterion the thresholding double ridge ratio criterion (TDRR). The detailed construction will be described in the following sections.

It is worth mentioning that selecting proper ridges is important for constructing a  good estimate. It would be argued that selecting ridges  be in spirit similar to selecting penalties in penalized criteria such as the BIC type criterion (Zhu et al., 2006). However, we can find  in the empirical studies that the criterion is much less sensitive to the ridge selection than the penalty selection in penalized criteria such as BIC and very much saves computational workload compared with BIC. More interestingly, when considering a sequence of local models with the structural dimension $q$, a consistent estimate $\hat q$ can be constructed with a wider range of ridges than that of penalties (see, e.g. Guo et al., 2015). Further, to better estimate the dimensionality, we give two ridges in the two rounds of ratios although in theory, the same ridge in the two rounds is also possible. To assist practical use, we recommend two values.

The materials of the paper are organized as follows.  Because estimating the target matrix and its eigenvalues is necessary, we first briefly review two popularly used methods: sliced inverse regression (SIR, Li, 1991) and discretization-expectation estimation (DEE, Zhu et al., 2010) in Section~2. In Section~3.1,  the thresholding double ridge ratio criterion (TDRR) is developed and its properties are described for sufficient dimension reduction with fixed and divergent number of covariates. The results about the local models are presented in Section~3.2. Section~4 exhibits its application to determining the number of common factors of  approximate factor models for ultra-high dimensional paradigms.  In Section~5, the numerical studies including simulations and real data examples are conducted to examine the performance of the new method and to compare with commonly used methods. Concluding remarks are included in Section~6 to discuss some possible research topics. Technical proofs are postponed to the appendix.

\section{ A Brief Review on Sufficient Dimension Reduction}
As described in the introduction, when we determine the structural dimension $q$ in sufficient dimension reduction (SDR), an estimation often requires the eigen-decomposition of a target matrix $M$ say.  We now briefly review two popular SDR methods: sliced inverse regression (SIR, Li, 1991) and discretization-expectation estimation (DEE, Zhu et al., 2010) although some recent developments require less technical conditions (see e.g. Li and Dong, 2009; Dong and Li, 2010; and  Guan, et al., 2016). Further for ease of exposition and notational convenience, we consider $Y$ to be scalar in the following sections. From DEE (Zhu et al. 2010), we can see that  it has no any difficulty to handle the $m$-dimensional response case.

\subsection{Sliced Inverse Regression}
Sliced inverse regression (SIR) method is an innovative idea
for obtaining  $\mathcal{S}_{Y|X}$ in regression. This promising method implicates the inverse regression that $X$ is regressed on $Y$.
Define the standardized covariates $X$ as the new covariates $Z=\Sigma^{-1/2}(X-u)$, where $u$ and $\Sigma$ denote the mean and the non-singular covariance matrix of $X$, respectively.
Under this linearity condition that the conditional mean $E(X|B^\top X)$ is a linear function of $B^{\top} X$ with the columns of $B \in \mathbb{R}^{p \times q}$ being any basis of the central space $\mathcal{S}_{Y|X}$, Cook (1998, Proposition 6.1) has justified that $\mathcal{S}_{Y|X}= \Sigma^{-1/2}\emph{S}_{Y|Z}$.
SIR is due to the fact that for any $y$, $ E(Z|Y=y) \subseteq \emph{S}_{Y|Z}$. Define the inverse mean subspace $\emph{S}_{E(Z|Y)} = \rm{Span}\{E(Z|Y=y), y \in \mathbb{R}\}.$
Then, under the coverage condition $\emph{S}_{E(Z|Y)} =\emph{S}_{Y|Z}$, the eigenvectors associated with non-zero eigenvalues of $\rm{Cov}\{E(Z|Y )\}$ make up a basis of  $\emph{S}_{Y|Z}$.

When an i.i.d sample  $\{x_i,y_i \}_{i=1}^n$ of $(X,Y)$ is available, from Li (1991), the SIR algorithm is  with the following steps.
Let $\bar{X}$ and  $\hat{\Sigma}$ be respectively  the sample mean and covariance matrix of $X$.
By partitioning the range of $Y$ into $H$ intervals, $I_1,\ldots, I_H$,  $\rm{Cov}\{E(Z|Y)\}$ is
approximately estimated by
\begin{eqnarray}
\widehat{\rm{Cov}}\{E(Z|Y)\}=\sum_{k=1}^H \hat{p}_k\bar{z}_k\bar{z}^\top_k,
\end{eqnarray}
where for $k=1,\ldots,H$, $p_k=P(Y\in I_k)$ and $z_k = E(Z|Y\in I_k)$ are  respectively  estimated by $\hat{p}_k= n_k / n$ with $n_k$ being the number of observations falling in the interval $I_k$ and the sample mean $\bar{z}_k$  in each interval $I_k$.

The spectral decomposition of $\widehat{\rm{Cov}}\{E(Z|Y)\}$, when $q$ is given, gets $q$ eigenvectors $\hat A_q =(\hat{\alpha}_1,\ldots, \hat{\alpha}_{q})$ associated with the $q$ largest eigenvalues. Hence, $\hat B_q = (\hat{\beta}_1, \ldots,\hat{\beta}_q) = \hat{\Sigma}^{-1/2} (\hat{\alpha}_1, \ldots, \hat{\alpha}_q)$ is a base matrix of the column space $\emph{S}_{Y|X}$ and $\hat{B}^{\top}X$ is the vector of the projected covariates.

\subsection{Discretization-Expectation Estimation}
To avoid the choice of the number $H$ of slices which sometimes has strong impact for the performance of the SDR estimation methods (see e.g., Zhu and Ng 1995; Li and Zhu 2007), Zhu et al. (2010) proposed the DEE method. The following are the steps of the SIR-based DEE algorithm.
\begin{enumerate}
\item [S1.] Dichotomize the response variable $Y$ into a set of binary variables by defining
 $Z(t)= I\{Y \leq t\}$, where the indicator function $I\{Y \leq t\}$ takes  value 1 if $Y \leq t$ and 0 otherwise.
\item [S2.] Let $\emph{S}_{Z(t)|X}$ denote the central subspace of $Z(t)|X$. When SIR is used, the related SIR matrix  $M(t)$ is  an $ p\times  p$ positive semi-definite matrix satisfying $\rm{Span}\{M(t)\} =\emph{S}_{Z(t)|X}$ in many cases.
\item [S3.] Let $\tilde Y$ be an independent copy of $Y$. The target matrix is $M=E\{M(\tilde Y)\}$. The matrix $ B$ consists of the eigenvectors associated with the nonzero eigenvalues of $M$.
\item [S4.] Obtain an estimate of $M$ as:
\begin{equation*}
M_{n}=\frac{1}{n}\sum^{n}_{i=1}M_n({y}_i),
\end{equation*}
where $M_n({y}_i)$ is an estimate of the SIR matrix $M({y}_i)$. When $q$ is given,  an estimate $ \hat B_q$ of $B$ consists of the eigenvectors  associated with the $q$ largest eigenvalues of $M_{n}$.
\end{enumerate}
It has been proved that $ \hat B_q$ is root-$n$
consistent to $B$. More details can be referred to Zhu et. al. (2010).
\section{Thresholding  Double Ridge Ratio Estimation}
\subsection{ The Motivation and  Criterion Construction}
We are now in the position to define the  criterion. Let $\lambda_1\ge \lambda_2 \ldots \ge \lambda_q> \lambda_{q+1}=\ldots=0$ be the eigenvalues of a target matrix $M$ at the population level. As we described before,  the first or the first two largest eigenvalues would be often dominating. Then we consider `standardizing' the eigenvalues by a strictly monotonic transformation: for $j=1,2,\ldots, p$,
$$
s_j:=\frac{{\lambda}_{j}}{{\lambda}_{j}+1}.
$$
All `standardized' eigenvalues are smaller than  $1$.  Then, $s_1\ge s_2\ge \ldots \ge s_q>s_{q+1}=\ldots =0$. Thus, the ratios
$$
s^*_j:=\frac{s_{j}^2}{s_{j+1}^2}-1
$$
have the following property: if defining $0/0$ as $1$,
\begin{eqnarray}\label{firstratio}
s^*_j=
\left\{\begin{array}{ll}
\left(\frac{s_{j}}{s_{j+1}}\right)^2 -1 \geq 0,
 & {\rm{\ for}}\  1\leq j < q,\\
+\infty, & {\rm{\ for}}\   j = q,\\
0/0-1= 0, & {\rm{\ for}}\  q+1 \leq j \leq p-1.
\end{array}\right.
\end{eqnarray}
In other words, the function $s^*_j$ about $j$ attains the maximum at $j=q$.
Based on this property, we could define a criterion to determine  $q$.

{\it Initial ridge ratio criterion:}
As is well known,  when $p<n$, all the  eigenvalues $\hat{\lambda}_{1} \geq \ldots \geq \hat{\lambda}_{p}$ of an estimated target matrix $M_n$ are usually non-zero, we then modify the above criterion such that at the sample level, the function can have the same property. Note that usually $\hat \lambda_j$ are convergent to the corresponding eigenvalues $\lambda_j$ at a certain rate that is typically  $O_p(\sqrt{p/ n})$ (Wu and Li, 2011). Thus, to avoid the instability of the ratios  $0/0$ for $j>q$, we add a ridge in the ratios ${s}^*_j$ for $j=1, \ldots p-1$ to define its sample version (see, Xia et al., 2015):
\begin{eqnarray}\label{tilderatio}
\hat{s}^*_j=\frac{\hat s^2_{j}+c_{1n}}{\hat s^2_{j+1}+c_{1n}}-1,
\end{eqnarray}
with $ \hat s_{j} = \frac{\hat{\lambda}_{j}}{\hat{\lambda}_{j}+1},$ for
$j=1,2,\ldots, p-1$ and the ridge  $c_{1n}$ is a value going to zero at a certain rate to be specified in the following theorem.  The choice of the ridge $c_{1n}$ is based on the  principle that $c_{1n}$ converges to zero at a slower rate than that of $\hat s^2_j \to 0$ for $j>q$. This can ensure that for $j=q$, $\hat s_j^*$ tends to $+\infty$ and for $j<q$, it is bounded from above by a constant. Therefore, for all $1\le j\le p-1$, we have that, at a certain convergence rate,
\begin{eqnarray}\label{firstratio1}
\lim_{n\to \infty}\frac 1{\hat s^*_j}=
\left\{\begin{array}{ll}
\frac 1{\left({s_{j}}/{s_{j+1}}\right)^2 -1} > 0,
 & {\rm{\ for}}\  1\leq j < q,\\
0, & {\rm{\ for}}\   j = q,\\
\ge 1, & {\rm{\ for}}\  q+1 \leq j \leq p-1.
\end{array}\right.
\end{eqnarray}
Two relevant references are Xia et al. (2015), and Zhu et al. (2016).
Thus,  theoretically the minimizer $\hat q$ of $1/\hat s_j^*$ over all $1\le j\le p-1$ can be proved to be equal to $q$ with a probability going to one. However,  the numerical studies that are not reported in this paper suggested that this criterion also underestimates the dimension $q$  as almost all existing methods performed. This is because the previous ratios before the $q$th could be smaller due to  the existence of the dominating eigenvalues. This is  particularly the case for local models (see Guo et al., 2015 and Zhu et al., 2016).  Thus, we define the second round of ridge ratios to further distinguish the value at $q$ from the others.

{\it Double ridge ratio criterion:} Consider the ratios between $s^*_j$ and $s^*_{j{+1}}$. Again  $0/0$ is defined as $1$ , and then
\begin{eqnarray}\label{secondratio}
\frac{s^*_{j+1}}{s^*_j}=
\left\{\begin{array}{ll}
C_j\geq 0,
 & {\rm{\ for}}\  1\leq j < q-1,\\
 +\infty, & {\rm{\ for}}\   j = q-1,\\
0, & {\rm{\ for}}\   j = q,\\
(0/0-1)/ (0/0-1)= 1, & {\rm{\ for}}\  q+1 \leq j \leq p-2.
\end{array}\right.
\end{eqnarray}
 Thus, $\frac{s^*_{q}}{s^*_{q-1}}$ is on peak and $\frac{s^*_{q+1}}{s^*_q}$ is at valley bottom. Unlike the ratios ${s}^*_{j}$ in the first round of ratios,  the  double ridge ratios could better separate the $q$th  ratio from the previous ratios $\frac{s^*_{j+1}}{s^*_{j}}$ for $1\le j\le q-2$. Again, to avoid the instability of the ratios $0/0$, we use the second ridge for the ratios between $\hat{s}^*_j$ and $\hat{s}^*_{j+1}$, to define
$({\hat{s}^*_{j+1}+c_{2n}})/({\hat{s}^*_j+c_{2n}}).$  When properly selecting the ridge $c_{2n}$, we expect that
\begin{eqnarray}\label{secondratio}
\frac{\hat{s}^*_{j+1}+c_{2n}}{\hat{s}^*_j+c_{2n}} \rightarrow  \frac{{s}^*_{j+1}}{{s}^*_j}.
\end{eqnarray}
Also, if following the classical methods, we could also use the minimizer of the ratios to define an estimate of $q$. Again, we do find that practically, the first or first two  ratios would still be local minima, and existing methods then underestimate $q$. Note that the true dimension could be the largest local minimizer and thus, we construct the following criterion to avoid the multiple local minima problem.

{\it Thresholding double ridge ratio criterion:}
 The maximum index among the indexes that make the ratios smaller than a threshold $0<\tau <1$ can be defined as an estimate. That is, the  dimension $q$ can be estimated by
\begin{eqnarray}\label{ratio}
\hat{q}=\arg\max_{1\leq j \leq p-2}\left\{j: \ \frac{\hat{s}^*_{j+1}+c_{2n}}{\hat{s}^*_j+c_{2n}}\leq \tau \right\}.
\end{eqnarray}

It is worth pointing out that this criterion is particularly  useful for local models in Section~3.2 below. This is because for local models, there will be more than one local minimum of the criterion and thus simply using minimizer or maximizer does not work.

To investigate the asymptotic properties, we assume the following regularity condition.
\begin{enumerate}
\item [A1.] When $p$ is fixed, $||M_n-M||=O_p(1/\sqrt{n})$,  where $||\cdot||$ is the Frobenius norm.
\end{enumerate}
This condition is very mild as almost all existing estimators have this root-$n$ consistency such as  SIR and DEE.

\begin{theorem}\label{theorem0}
In addition to  Condition~A1,  assume that $c_{1n} \rightarrow 0$, $c_{2n} \rightarrow 0$,  $c_{1n}c_{2n}n \rightarrow \infty$, $0<\tau<1$ and $p$ is fixed. Then  the estimate $\hat{q}$ by (\ref{ratio}) is equal to $q$ with a probability going to $1$.
\end{theorem}
How to choose an optimal parameter $\tau$ plays an important role in practice although it is not an issue in theory.  It is understandable that  if $\tau$ is too close to $0$, the method tends to underestimate  $q$ and if $\tau$ is too close to $1$, to overestimate $q$. By the rule of thumb, we  recommend $\tau=0.5$. The numerical studies in the late section  support this choice.

When $p$ is divergent as the sample size $n$ goes to infinity, the root-$n$ consistency of $M_n$ is usually no longer true. Therefore, we need to adjust the ridges. To this end, we first check the convergence rate of the estimated matrix $M_n$ to the target matrix $M$.
Assume the following regularity condition:
\begin{enumerate}
\item [A2.] $||M_n-M||=O_p(\sqrt{p/n})$, where $||\cdot||$ is the Frobenius norm.
\end{enumerate}
This rate is also common for the SDR estimates in the literature, see Wu and Li (2011). 

\begin{theorem}\label{theorem2.2}
Assume that  Condition A2 holds, $p=O_p(n^{\alpha})$ with $0\leq\alpha<1$, $q$ is fixed, $c_{1n} \rightarrow 0$, $c_{2n}n \rightarrow 0$  and  $c_{1n}c_{2n}n/p \rightarrow \infty$. Then the estimate $\hat{q}$ by (\ref{ratio}) is equal to $q$ with a probability going to $1$.
\end{theorem}

\subsection{Estimating $q$ in A Sequence of Local Models }
Recall in Section~1, we consider the sequence of multiple-index models in (\ref{sequence}) as:
\begin{eqnarray}
 Y_n=G(B^{\top}_1X)+ C_n g(B^{\top}X)+\varepsilon.
\end{eqnarray}
Rewrite
 $S_{E(Y|X)}=\rm{Span}(B_1)$ and $S_{E(Y_n|X)}=\rm{Span}(B_1, B)$, respectively. 
From Guo et al. (2015), we can see that identification of  $\rm{Span}(B_1, B)$ plays an important role for the omnibus property   of a test. It is however a very challenging problem. In a spacial setting that $B_1$ is a vector ($q_1=1$) contained in $\rm{Span}(B)$,  Guo et al. (2015) used the BIC type criteria (Wang and Yin, 2008; Zhu et al., 2010) to estimate the dimension $q$ of $\rm{Span}(B_1, B)$. They proved that when $C_n=O(n^{-1/4}h^{-1/2})$ where $h \to 0$ is the bandwidth in the kernel estimate of nonparametric regression function,  the conclusion that $\hat q = q$ with a probability going to $1$  no longer holds.  When the new method is used, we now give a systematic investigation to examine that at what rate of $C_n\to 0$, $\hat q$ is still consistent and at what rate it fails to have this property.  The results are stated in the following theorem.

\begin{theorem}\label{theorem-seq}
Under Conditions~A3--B4 in the Appendix and the sequence of local models in $(\ref{sequence})$ with fixed $p$,
 the estimate $\hat{q}$ by (\ref{ratio}) has the following consistency.
\begin{itemize}
\item [(I)]  If $C_n =O(n^{-1/2})$, $c_{1n} \rightarrow 0$, $c_{2n} \rightarrow 0$,
and $c_{1n}c_{2n}n \rightarrow \infty$, $P(\hat{q}= q_1)\rightarrow 1$.
\item [(II)] If $C_n \approx n^{-\alpha} $ with $0<\alpha<1/2$, when   $c_{1n}=o_p(C^2_n)$, $c_{2n} \rightarrow 0$, $c_{1n}c_{2n}n \rightarrow \infty$ and $c_{1n}c_{2n}/C^4_n \rightarrow \infty$, $P(\hat{q}= q)\rightarrow 1$.
\end{itemize}
\end{theorem}

The results give a full picture about the possible underestimation and also  show that the new criterion has a significant improvement over the BIC criteria used in Guo et al. (2015). That is, when $C_n \approx n^{-\alpha}$ with $0<\alpha <1/2$, the new method could well identify the true dimension $q$ of the local models if the ridges could be chosen properly. In contrast, whenever the rate is of the order  $n^{-1/2}$ or faster, underestimation  always occurs.

\underline{}\section{ Determining The Number of Common Factors in Approximate Factor Models }
%
Factor analysis is a useful statistical tool because it takes dimension reduction into consideration, see Bai and Ng (2002). An important issue is to determine the number of common factors in high-dimensional factor models. Consider the ultra-high dimensional approximate factor model as (Fan et al., 2008; Wang, 2012):
\begin{eqnarray}\label{factor0}
Y_{i}=\mathbf{B}F_i+U_{i}, \ \ i=1,2,\ldots,n,
\end{eqnarray}
where $Y_i$ is $p$-dimensional observation, $F_i=(f_{i1},\ldots, f_{id})^{\top}\in \mathbb{R}^d$  is an unknown $d$-dimensional column vector of common factors, $\mathbf{B}=(B_1,\ldots, B_p)^{\top}\in \mathbb{R}^{p\times d}$ with $B_j$ being the  loading vector for the $j$th component of $Y_i$ and
$U_{i}=(u_{i1},u_{i2},\ldots, u_{ip})^{\top}$ is the idiosyncratic factor vector that is uncorrelated with the common factors. The model (\ref{factor0}) can be rewritten in matrix form as:
\begin{eqnarray}\label{factor1}
\mathbf{Y}=\mathbf{B}\mathbf{F}+\mathbf{U},
\end{eqnarray}
where  $\mathbf{Y}=(Y_{1},\ldots, Y_{n})$ is the $p\times n$ observation matrix, $\mathbf{B}$ is the $p\times d$ factor loading matrix,  $\mathbf{F}=(F_1,\ldots, F_n)$ is the $d\times n$ common factor matrix and $\mathbf{U}=(U_1,\ldots, U_n)$  is the $p\times n$  idiosyncratic factor matrix.

Let $m=\min\{n,p\}$ and $\hat{\lambda}_{1} \geq \ldots \geq \hat{\lambda}_{p}$ stand for the eigenvalues of the matrix $\mathbf{Y}\mathbf{Y}^{\top}/(np)$.
As documented by Wang (2012) and Ahn and Horenstein (2013), the first $d$ eigenvalues are $O_p(1)$ and the rest are at most $O_p(1/m)$, we have a similar definition as that in (\ref{tilderatio}):
\begin{eqnarray}\label{tilderatio-factor}
\hat{s}^*_j=\frac{\hat s_{j}+c_{n,p}}{\hat s_{j+1}+c_{n,p}}-1\  {\rm{with}}\ \hat s_{j} = \frac{\hat{\lambda}_{j}}{\hat{\lambda}_{j}+1},
 \  \  j=1,2,\ldots, p-1,
\end{eqnarray}
where the choice of $c_{n,p}$ is discussed in the following theorem.
Comparing the definition (\ref{tilderatio-factor})  with (\ref{tilderatio}), we use $\hat s_{j}$ rather than $\hat s^2_{j}$ because $\hat s_{j}$ for $j>d$ converges to zero at a convergence rate of second order $O_p(1/m)$ rather than  the rate of first order $O_p(1/\sqrt{m})$.
Thus, the number $d$ of common factors can be estimated by
\begin{eqnarray}\label{ratio-factor}
\hat{d}= \arg\max_{1\leq j \leq p-2}\left\{j:\    \frac{\hat s^*_{j+1}+\tilde{c}_{n,p}}{\hat s^*_j+\tilde{c}_{n,p}}\leq \tau \right\},
\end{eqnarray}
where $0<\tau<1$ is some constant. As  aforementioned in Section~3,
the threshold value $\tau=0.5$ is adopted. The following  theorem shows that the estimate is consistent even for ultra-high dimension $p$.
\begin{theorem}\label{theorem-factor}
Under Conditions~B1--B3 in the Appendix, $c_{n,p}\rightarrow 0$, $\tilde{c}_{n,p} \rightarrow 0$, and $c_{n,p}\tilde{c}_{n,p}m \rightarrow \infty $, the estimate $\hat{d}$ defined by (\ref{ratio-factor}) satisfies that $P(\hat d=d)\to 1$ as $n \to \infty$.
\end{theorem}

\underline{}\section{Numerical Studies}
In this section, the simulation studies are carried out to examine the performance of the method and to compare with some popularly used competitors in the literature.

\subsection{Simulations for Structural Dimension Determination}
In the simulations, every experiment is repeated  $500$ times and the sample sizes are 200, 400 and 800 respectively. The competitors include the ridge-type ratio estimation (RRE, Xia et al., 2015) that is defined as:
\begin{eqnarray}
\hat{q}=\arg\min_{1\leq j \leq p-1}\left\{ \frac{\hat{\lambda}_{j+1}+c_n}{\hat{\lambda}_j+c_n}\right\};
\end{eqnarray}
the BIC-type estimation (Zhu et al., 2006 and Zhu et al., 2010) that has the form:
 \begin{eqnarray*}
\hat{q}=\arg\min_{1\leq j \leq p}
\left\{ \frac{n\sum_{l=1}^j\{\log(1+\hat{\lambda}_l)+\hat{\lambda}_l\}}
  {2\sum_{l=1}^p\{\log(1+\hat{\lambda}_l)+\hat{\lambda}_l\}}-\alpha_n\frac{j(j+1)}{p}\right\},
\end{eqnarray*}
where  $\alpha_n$ is the penalty term;
and the ratio estimation (RE, Luo et al., 2009) that  is defined as:
\begin{eqnarray*}
\hat{q}=\arg\min_{1\leq j \leq d_{max}}\left\{ \frac{\hat{\lambda}_{j+1}}{\hat{\lambda}_j}\right\},
\end{eqnarray*}
where $d_{max}$ is a  predetermined maximum number.  For easy comparison, the recommended values in the respective methods are used. Specifically, for RRE (Xia et al., 2015), the ridge value $c_n=\log(n)/(10\sqrt{n})$; for RE (Luo et al., 2009), the maximal number $d_{\max}=10$; for the BIC type criterion (Zhu et al., 2010), the penalty term $\alpha_n=\sqrt{n}$; for our TDRR, two ridge values  $c_{1n}=\log{n}/(10\sqrt{n})$ and $c_{2n}=\log{n}/(5\sqrt{n})$  are recommended. We also make a comparison with the sequential testing method (ST) (see Li, 1991; Bura and Cook, 2001).

\textbf{Example 1.} Data are generated from the following model:
\begin{eqnarray*}
Y=X_1/\{0.5+(X_2+1.5)^{1.5}\}+ X^3_3/4+\varepsilon;
\end{eqnarray*}
where $\varepsilon\sim N(0,\sigma^2)$ with $\sigma=0.2$, $X\sim N(\textbf{0},\textbf{I}_p)$, which is independent of $\varepsilon$. In this example, $q=3$, the dimension $p$ varies from 5 to 40. In this example, we adapt the SIR-based DEE procedure to obtain the estimated target matrix $M_{n}$. The results of the four methods under different combinations of sample sizes and dimensions of covariates $X$ are presented in Table ~\ref{model1-1}.

\begin{center}
Table ~\ref{model1-1} about here 
\end{center}

The results in the tables suggest that TDRR performs uniformly the best among all the competitors. When the sample size $n$ becomes larger, the proportion of correct decisions for our method increases higher reasonably. It is clearly observed that  no matter how large the dimension $p$ is,  RE and RRE estimate $q$ to be $1$ even for the large sample size $n = 800$.  BIC performs much better than  RE and RRE in high dimensional scenarios (especially when $p=40$). An interesting observation   is that the dimensionality is blessing BIC. Its performance gets better with higher dimension. But even though, it still  very much underestimates $q$ in the low dimensional cases. This means that the BIC method is not robust against different dimensions. As showed in Zhu et al. (2006), it is not easy to select proper penalty value $\alpha_n$. Overall, the performance of TDRR is much better. We also note 
that when the sample size is small ($n=200$), TDRR estimates $q$ to be $2$ with fairly high proportions showing that TDRR is also underestimation as the true structural dimension $q=3$, but  the proportions of $\hat{q}\ge 2$ is  always more than $80\%$ in all scenarios.

\begin{figure}
  \centering
  \includegraphics[width=7.250cm]{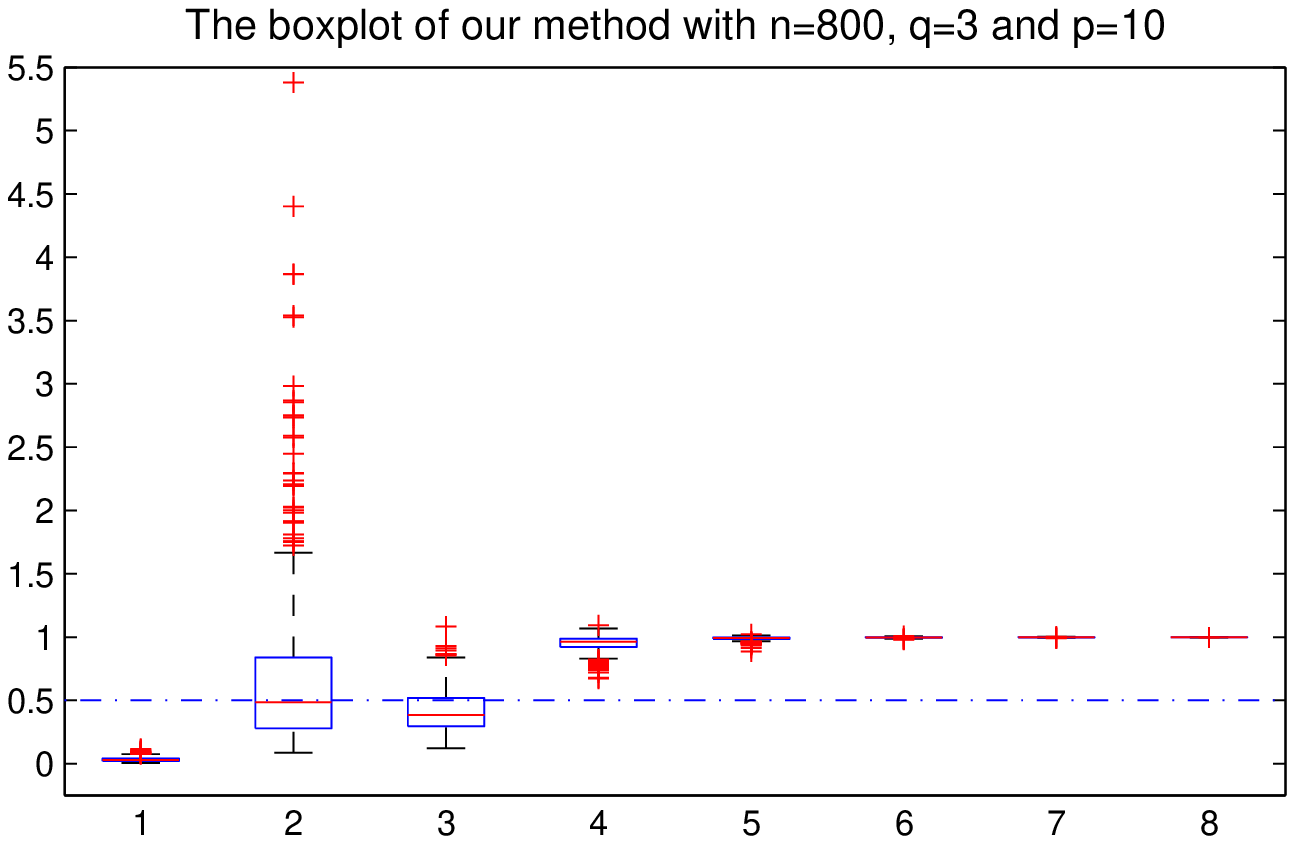}
  \includegraphics[width=7.250cm]{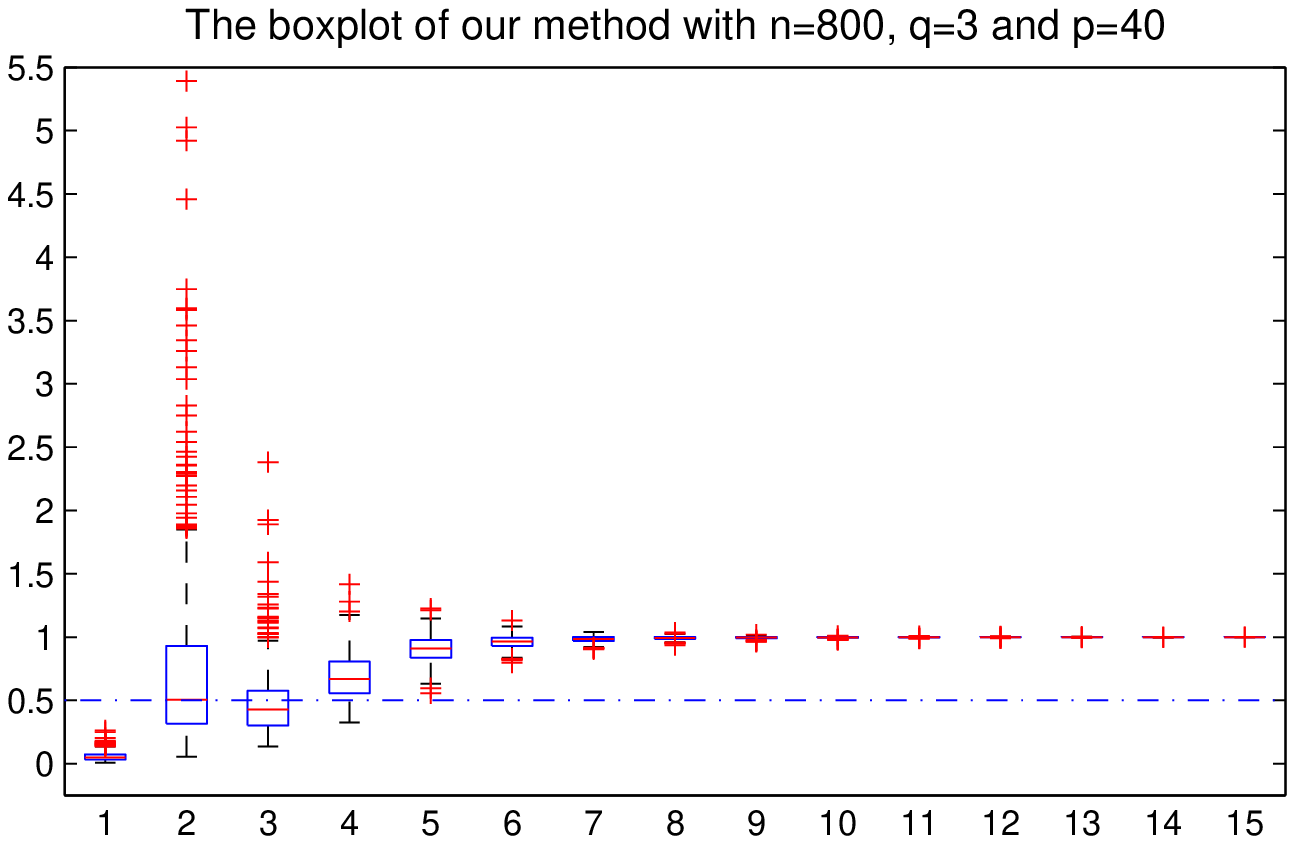}
  \includegraphics[width=7.250cm]{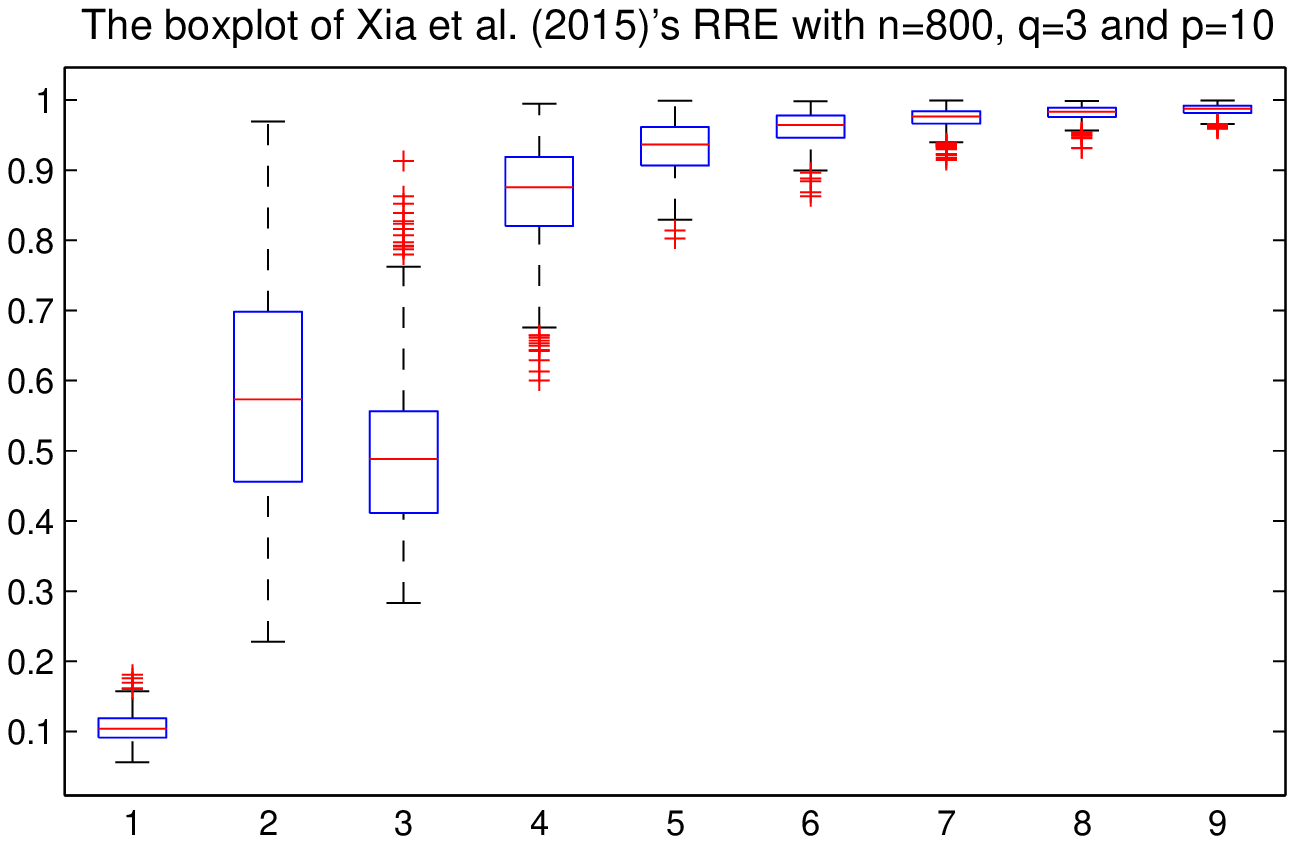}
  \includegraphics[width=7.250cm]{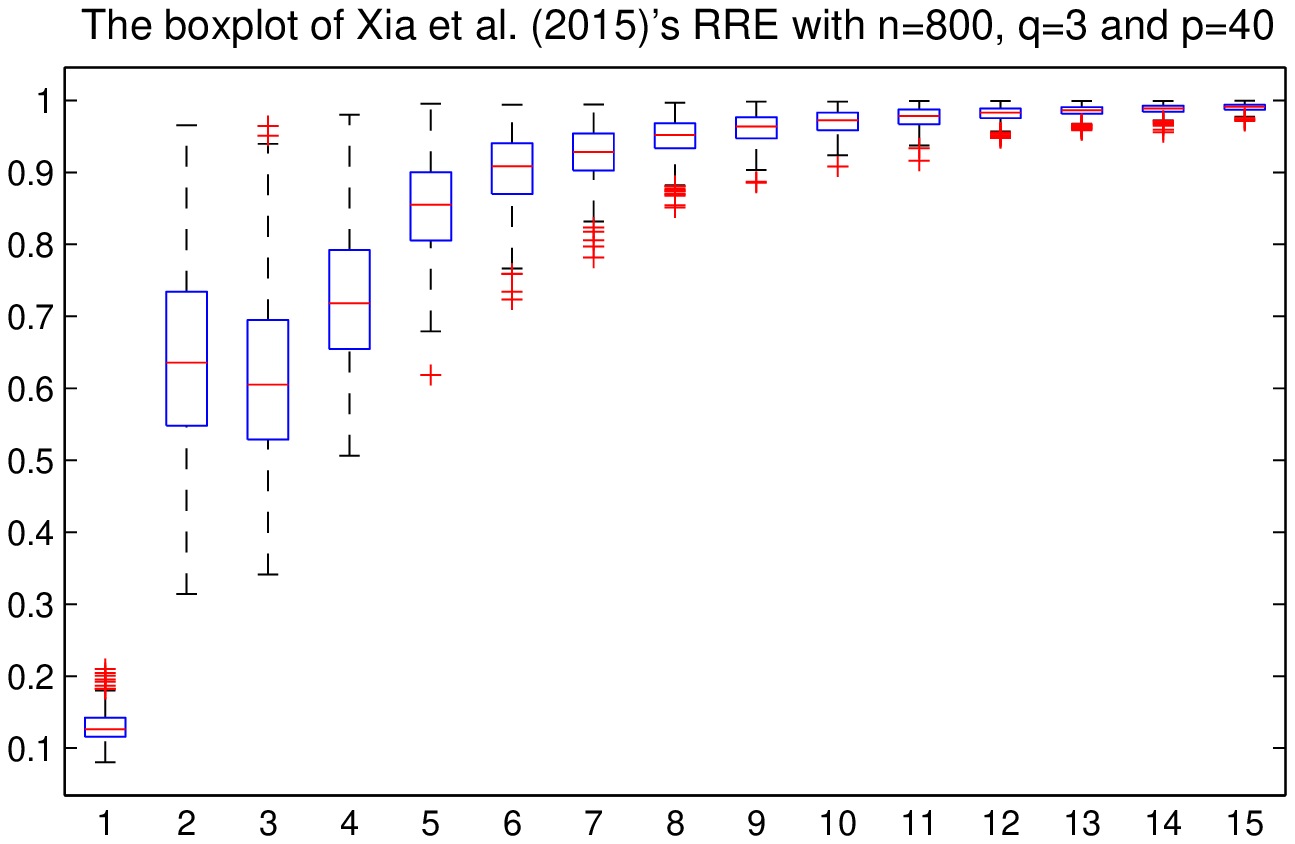}
  \includegraphics[width=7.250cm]{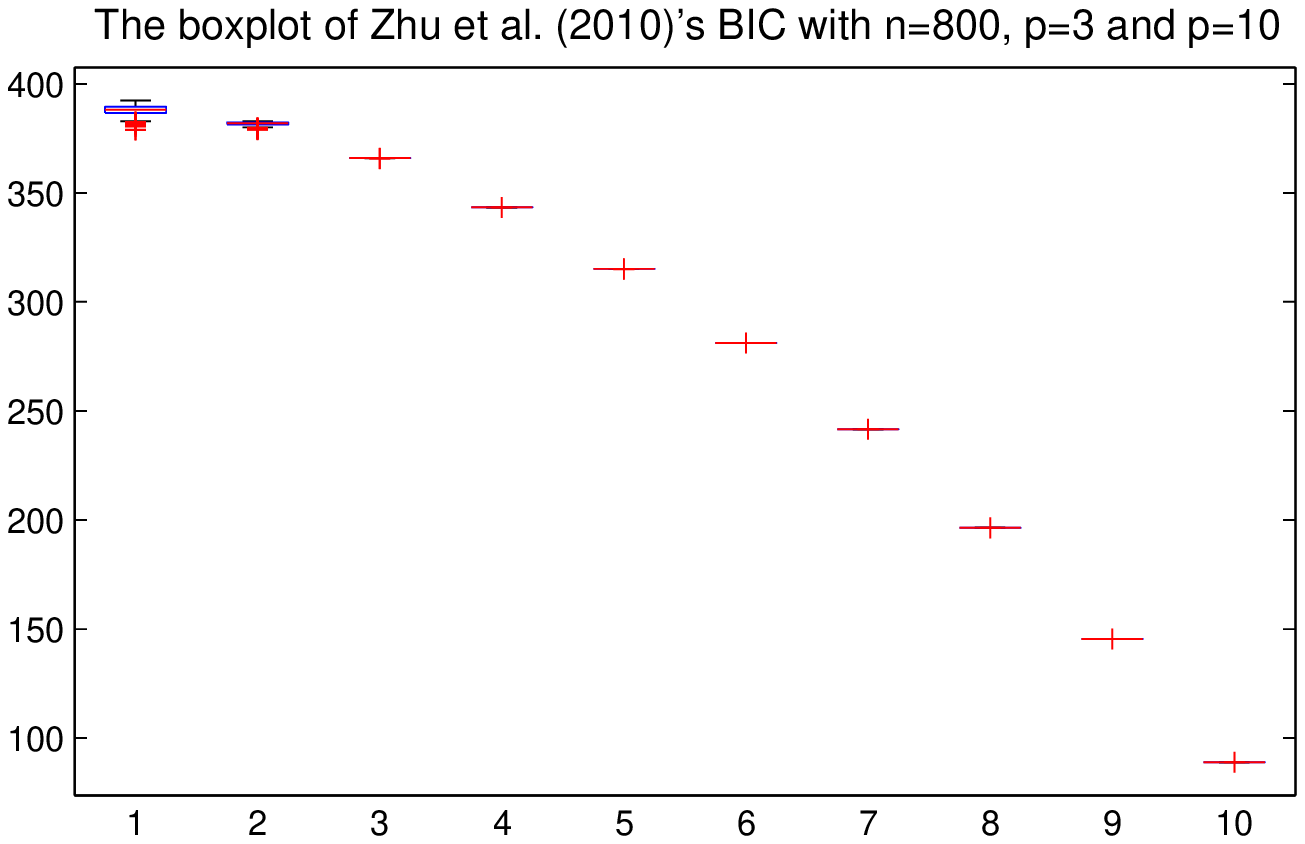}
  \includegraphics[width=7.250cm]{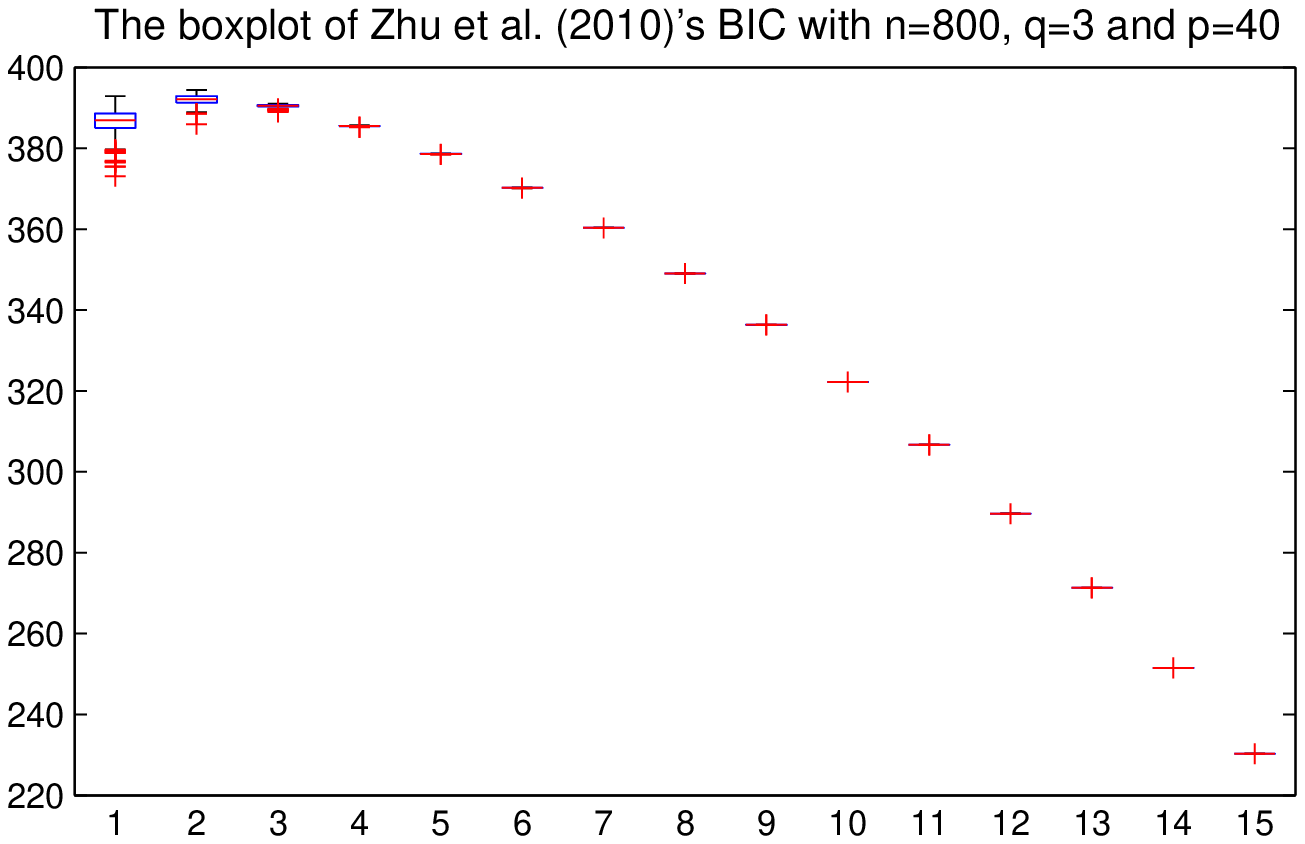}
  \includegraphics[width=7.250cm]{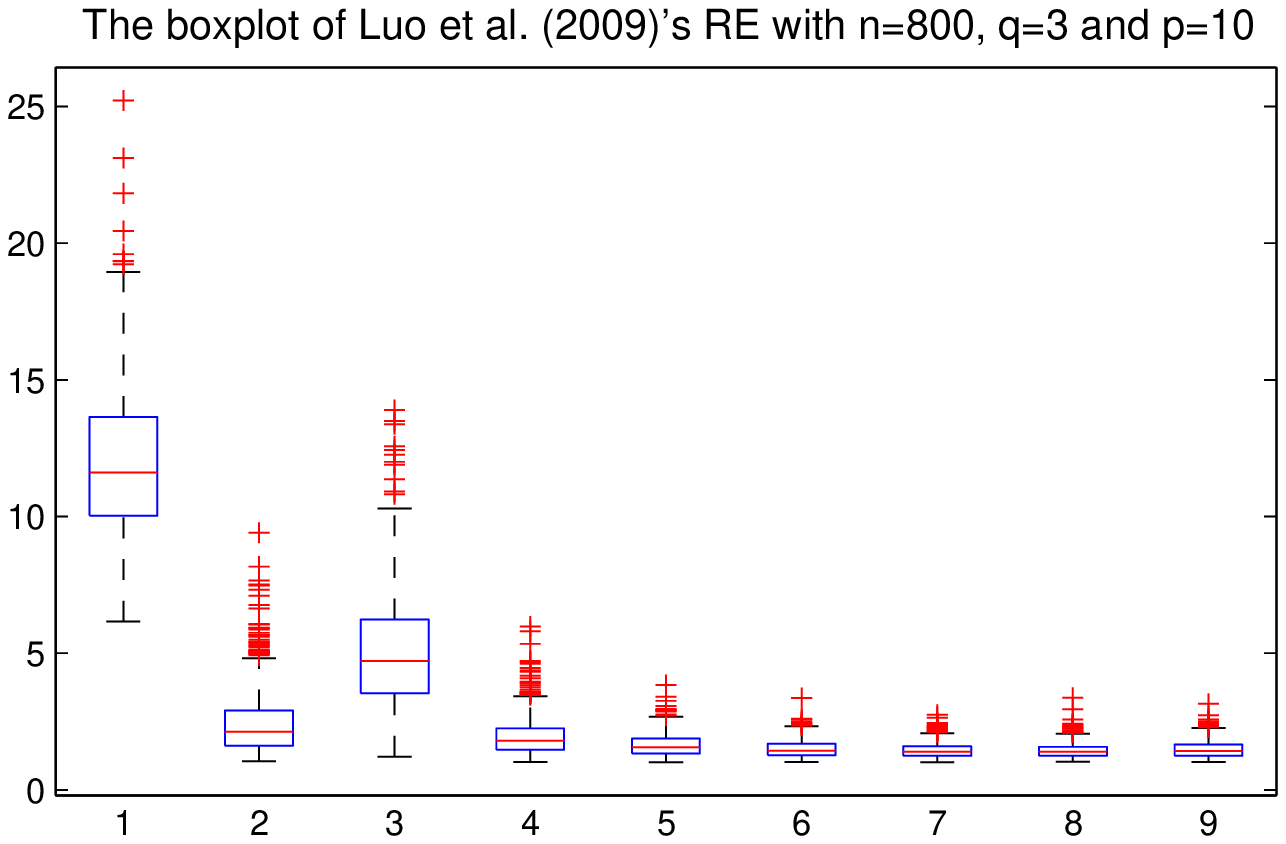}
  \includegraphics[width=7.250cm]{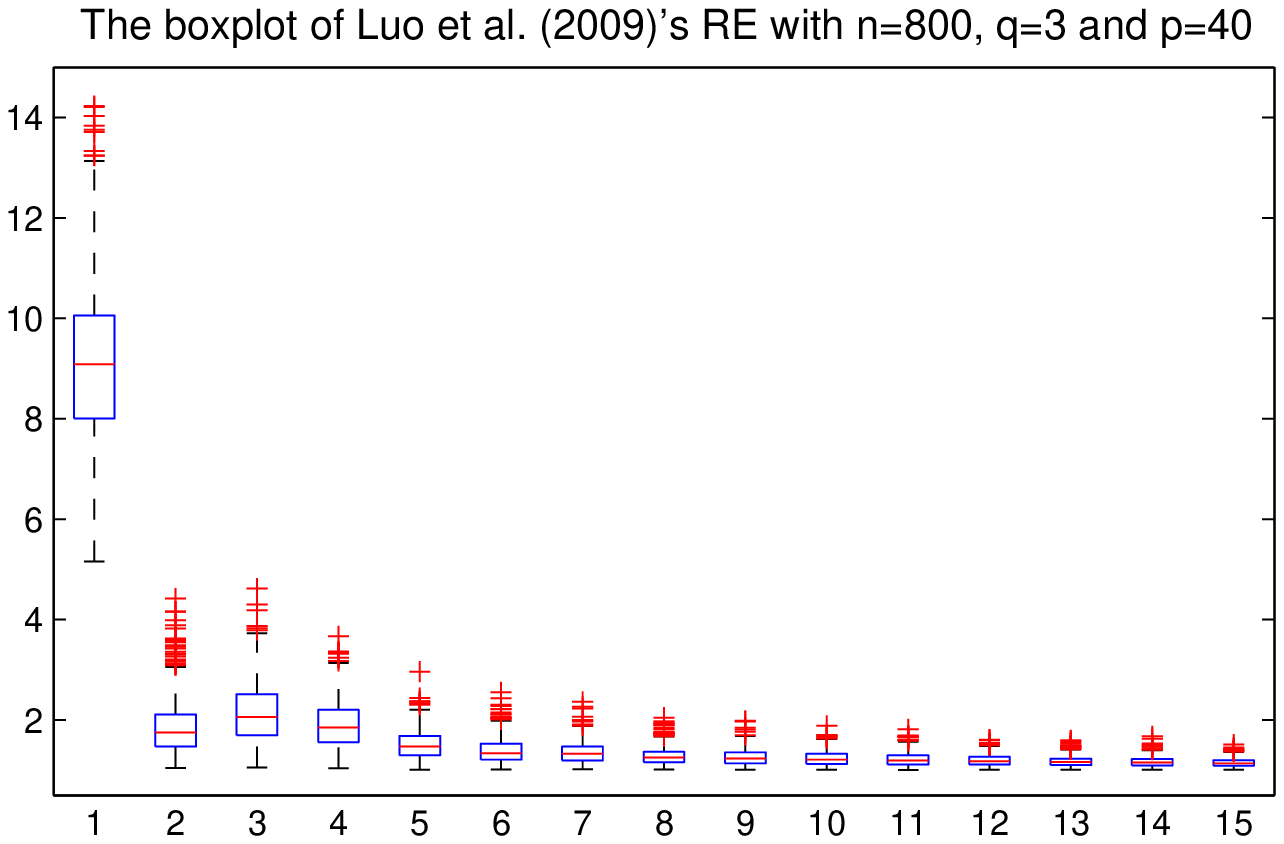}\\
  \caption{Boxplots  of the first 15 components of TDRR $\{\hat{s}^*_{j+1}+\frac{\log{n}}{5\sqrt{n}}\}/\{\hat{s}^*_j+\frac{\log{n}}{5\sqrt{n}}\}$ (TDRR),  $\{\hat{\lambda}_j+\frac{\log{n}}{10\sqrt{n}}\}/\{\hat{\lambda}_{j+1}+\frac{\log{n}}{10\sqrt{n}}\}$ (RRE), $\frac{n\sum_{l=1}^j\{\log(1+\hat{\lambda_l})+\hat{\lambda}_l\}}
  {2\sum_{l=1}^p\{\log(1+\hat{\lambda_l})+\hat{\lambda}_l\}}-\sqrt{n}\frac{j(j+1)}{p}$ (BIC) and
   $\frac{\hat{\lambda}_j}{\hat{\lambda}_{j+1}}$ (RE), for $j=1,\ldots, 15$ in Example~1. }\label{figure-reduction}
\end{figure}

To further reveal the reasons for the above phenomena, we draw the boxplots of  the first fifteen components of the four criteria in Figure~\ref{figure-reduction} when $n=800$. We have the following findings. First, the largest eigenvalue is shown to be dominating such that the first ratios of the four criteria attain either minimum or maximum. TDRR also has this problem.
 Thus, when using minimizer or maximizer to define an estimate of $q$, it should be $1$ rather than the true value $q=3$ and thus the underestimation is inevitable. From the first row of the plot for TDRR, we can see clearly the advantage of thresholding. As the third ratio attains a local minimum and thus, the maximum index of the ratios under the control by the thresholding value $\tau$ can be $\hat q=3$. This reveals the reason why TDRR has a much  better performance than the others.

 \begin{center}
Figure~\ref{figure-reduction} about here
\end{center}

\textbf{Example 2.}  Consider the  model used in Zhu et al. (2006) so that we can compare with the BIC method and sequential testing method (ST) (see Li, 1991; Bura and Cook, 2001) in a fair way:
\begin{eqnarray*}
Y=X_1\times (X_2+X_3+1)+ \varepsilon,
\end{eqnarray*}
where $\varepsilon\sim N(0,\sigma^2)$ with $\sigma=0.5$, $n=400,\ 800$ and $X \sim N(\textbf{0},\textbf{I}_p)$ which is independent of  $\varepsilon$. In this model, $q=2$. The dimension $p$ of covariates varies from $10$ to $40$ in the simulation.

As done in Zhu et al. (2006), we also use the SIR matrix (Li, 1991) as the estimated target matrix $M_n$.   Here the number of slices is selected to be $10$ as many references used (e.g. Li, 1991; Zhu et al., 2006).   The results are reported in Table~\ref{Example2}.
\begin{center}
Table~\ref{Example2} about here
\end{center}
In this table, the results of BIC and ST are excerpted from Zhu et al. (2006).
We can see that for the three methods, the proportions of correct decisions reasonably decrease with the increasing dimension $p$. When the sample size is large ($n=800$),  TDRR  significantly better performs than the two competitors, particularly when $p=40$. ST performs the worst, particularly for large $p$ ($p\geq 30$) with less than 35\% proportions of correct decisions. This implies that ST can not work  well in  high dimensional scenarios.

To check the performance of the proposed method for  local models and compare with the other methods, we use the following Examples~3 and~4.

\textbf{Example 3.} Consider the following sequence of models:
\begin{eqnarray*}
Y=(X_1+X_2)+ a*(X_1+X_2)X^{1.5}_3+\varepsilon,
\end{eqnarray*}
 where $\varepsilon\sim N(0,1)$ and  $X \sim N(\textbf{0},\textbf{I}_p)$  are mutually independent. the value $a= 2/n^{1/4}$ responds to a local model with  $q=2$ where the directions are $\beta_1=(1, 1, 0, \ldots, 0)^{\top}$ and $\beta_2=(0, 0, 1, 0, \ldots, 0)^{\top}$. When $a=0$, the structural dimension of the hypothetical model (its limiting model) is $q_1=1$ where $\beta_1=(1, 1, 0, \ldots, 0)^{\top}$.

\textbf{Example 4.}  We also conduct the simulation with  a local model which was designed by Guo et al. (2015) in studying the power performance of their test:
\begin{eqnarray*}
Y=0.25\exp(2X_1)+ aX^3_2+\varepsilon,
\end{eqnarray*}
where the structural dimension $q=2$ with $\beta_1=(1, 0, 0, \ldots, 0)^{\top}$ and $\beta_2=(0, 0, 1, 0, \ldots, 0)^{\top}$, $X \sim N(\textbf{0},\textbf{I}_p)$ and $\varepsilon\sim N(0,1)$ are mutually  independent, and $a=1/n^{1/4}$.  When $a=0$, the hypothetical model is with the structural dimension $q_1=1$ where $\beta_1=(1, 0, 0, \ldots, 0)^{\top}$. The results are reported in Table~\ref{Example4}.

\begin{center}
Table~\ref{Example4} 
about here
\end{center}

The following conclusions  can be arrived from the simulation results  in Table~\ref{Example4}. First, RE works better than RRE and BIC. But  the three competitors RE, RRE and BIC  very much underestimate the structural dimension to be $1$ even in large sample size cases ($n=800$). This again confirms the underdetermination problem for the competitors.  TDRR  performs much better than the other competitors.  

Additionally,  we have the following findings through the comparisons between Table~\ref{Example4} with Table~\ref{model1-1}. 
First,  for larger structural dimension  TDRR is slightly more difficult to identify, but the performance of the other competitors are  seriously deteriorated. Second, the large sample size $n$ can reasonably improve the performance of TDRR, but the competitors still seriously underestimate the structural dimension even with $n=800$. RRE and RE tend to estimate $q$ to be $1$ rather than the true value $3$. Third, although for TDRR, the proportion of correct decisions about $q$  decreases as the dimension $p$ increases, the frequencies that the estimated dimension $\hat q$ is more than the true value $3$ is more than 75\% when $n=800$. Fourth, TDRR does not seriously overestmiate $q$ such that the working models are still parsimonious enough. Again, the dimensionality blessing for BIC  exists.  This interesting phenomenon deserves a further study although BIC would  not  be recommendable when the dimension is very high due to its computational complexity.

Overall, TDRR is an efficient method to determine the structural dimension. Also, it is robust to different models, while the other competitors are not.

\subsection{Simulations for Approximate Factor Model}
Now we turn to investigate the performances of TDRR for approximate factor models.
Similarly, we consider the three aforementioned competitors:  RRE (Xia et al., 2015), RE (Wang, 2012) and BIC (Zhu et al., 2010) that was also used in comparison in Xia et al. (2015).
The ridge values $c_{n,p}= \log(m)/(10\sqrt{m})$ and $\tilde{c}_{n,p}=\log(m)/(5\sqrt{m})$ with $m=\min\{p, n\}$ are selected in our method and the maximum number $d_{max}=\min\{p, n\}/2$ is used to make  RE  more efficient in estimation. Following the suggestions of Xia et al. (2015), we respectively choose the ridge value $c_{n}=\log(n)/(10n)$ in RRE and the penalty term $C_n=\log(n)$ in BIC.

\textbf{Example 5.} The model is adopted from Fan et al. (2015) but with a more sophisticated
factor structure. The model is as follows:
\begin{eqnarray*}
\mathbf{Y}_{t}=\mathbf{B}\mathbf{F}_t+\mathbf{U}_{t}, \ t=1,2,\ldots,n,
\end{eqnarray*}
where $\mathbf{U}_t$ is from the standard multivariate normal distribution $N(0, \textbf{I}_p)$
and $\mathbf{F}_t$ is from multivariate $t$-distribution $t(\nu, \Sigma_{i})$ with degrees of freedom $\nu$ and covariance matrix $\Sigma_{i}$ for $i=1,2$ and $3$, where
$\Sigma_{1}=\textbf{I}_d$, $\Sigma_2=(\sigma_{ij}^{(2)})_{d\times d}$ and
$\Sigma_3=(\sigma_{ij}^{(3)})_{d\times d}$ with the elements respectively
\begin{eqnarray*}
\sigma_{ij}^{(2)}=I(i=j)+ 0.8^{|i-j|} I(i \neq j) \ \
{\rm and} \ \ \sigma_{ij}^{(3)}=I(i=j)+ 0.7I(i \neq j).
\end{eqnarray*}
Every row of $\mathbf{B}$ is from $N(0, \textbf{I}_d)$ with the number of factors $d=4$ in this example. We vary $p$ with the change of the sample size $n$: $p=n, 2n$ and $4n$.
Smaller $\nu$ associates with  heavier tail and $\nu=\infty$ makes normal tail as it corresponds to a multivariate normal distribution $N(0, \Sigma_i)$.
Two different degrees of freedom $ \nu = 2.5$ and $\nu=\infty$ are chosen, representing respectively heavy tail and normal scenarios.

\begin{center}
Tables~\ref{factor-N} and~\ref{factor-T} about here
\end{center}

\begin{figure}
  \centering
  \includegraphics[width=7.250cm]{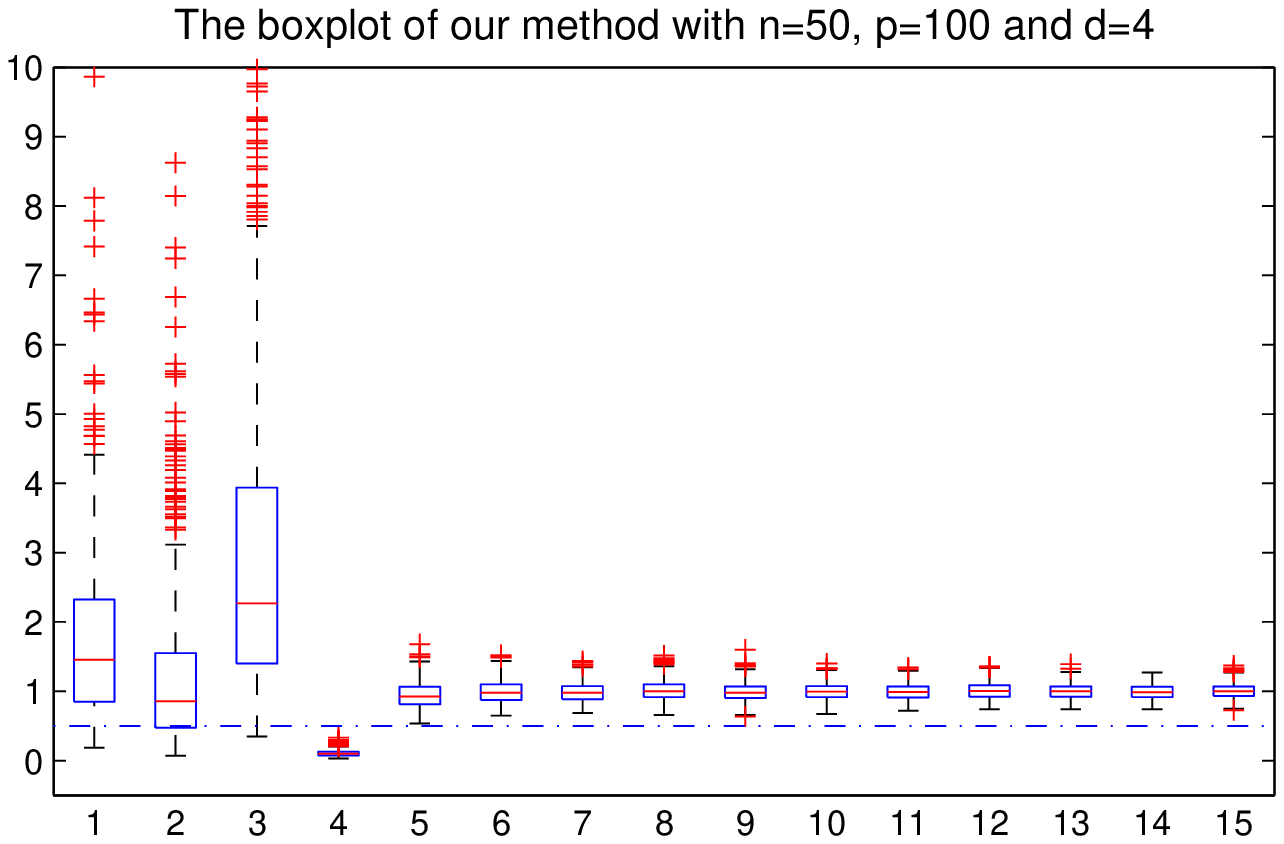}
  \includegraphics[width=7.250cm]{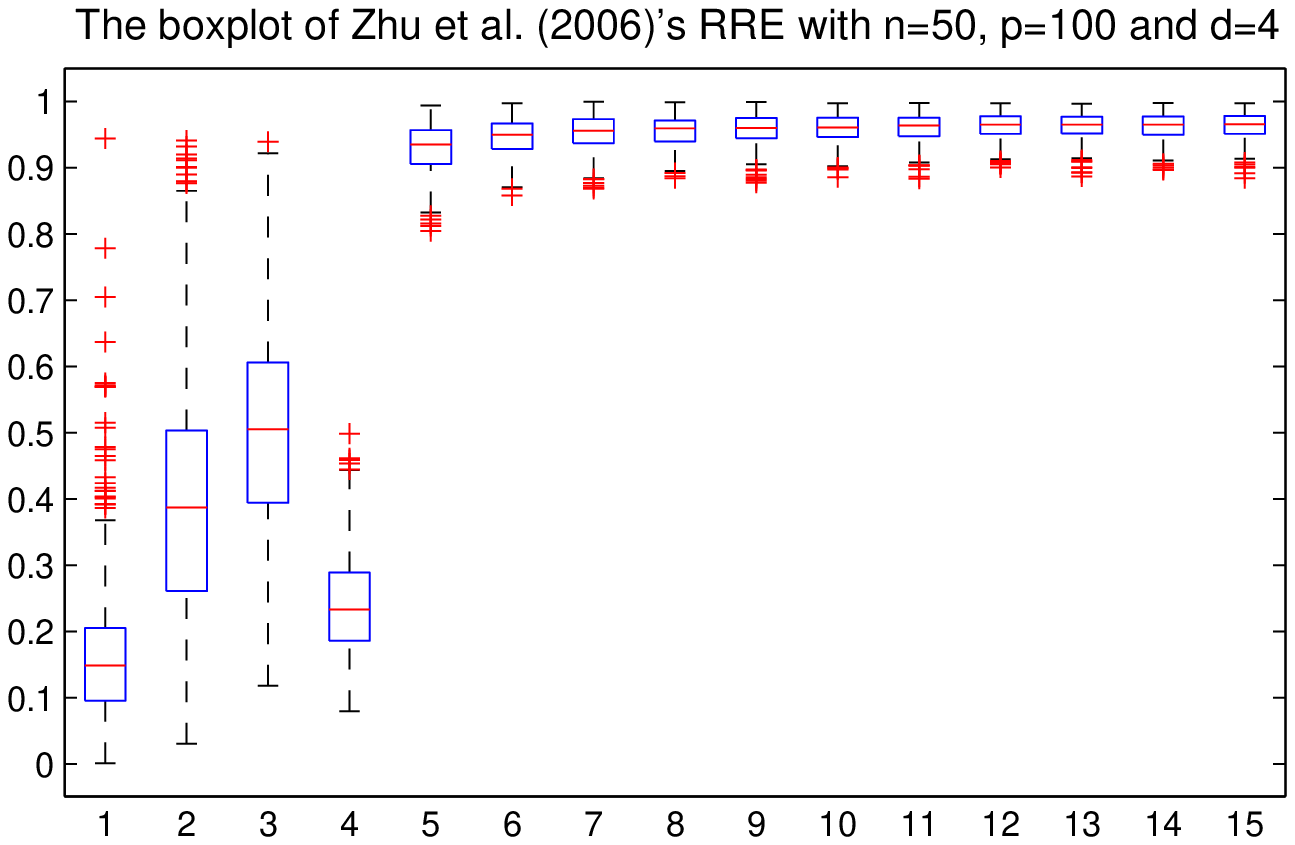}
  \includegraphics[width=7.250cm]{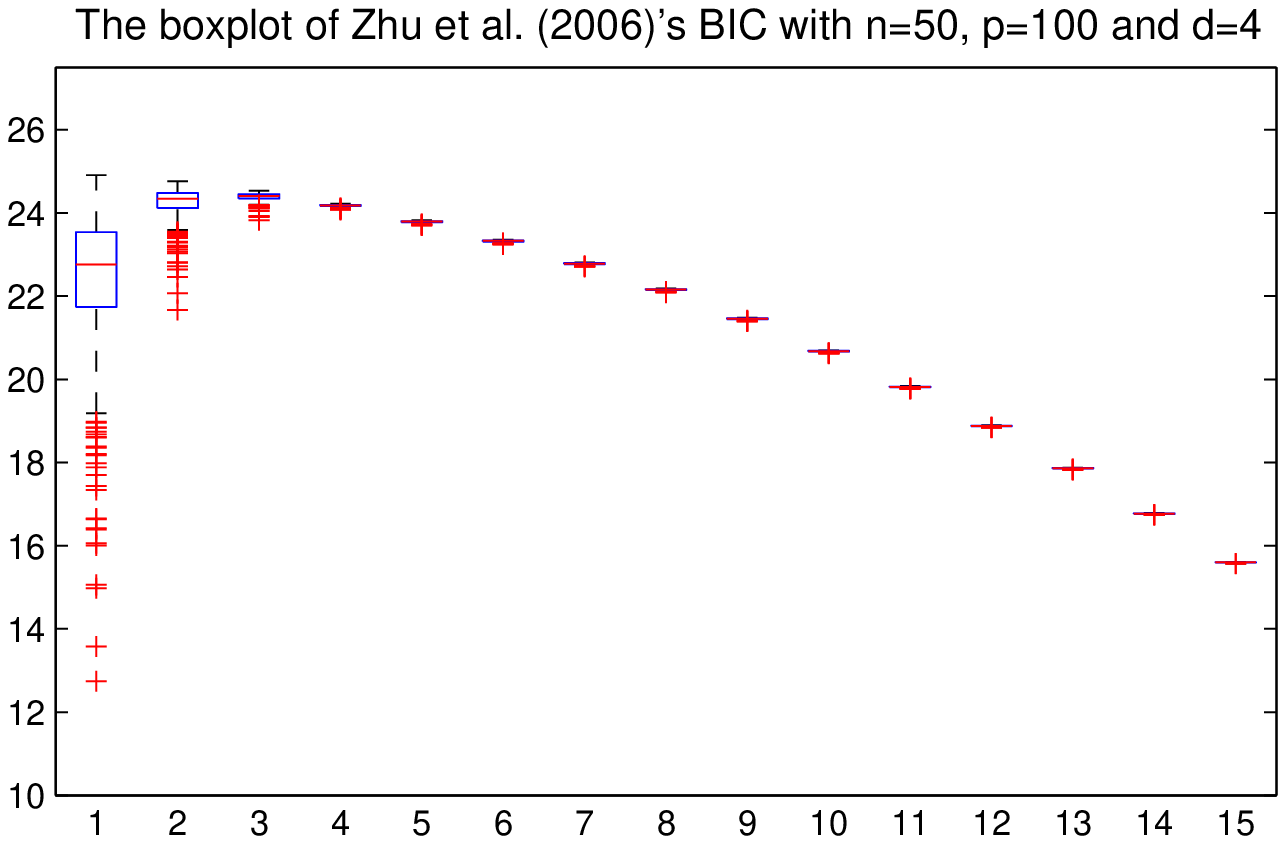}
  \includegraphics[width=7.250cm]{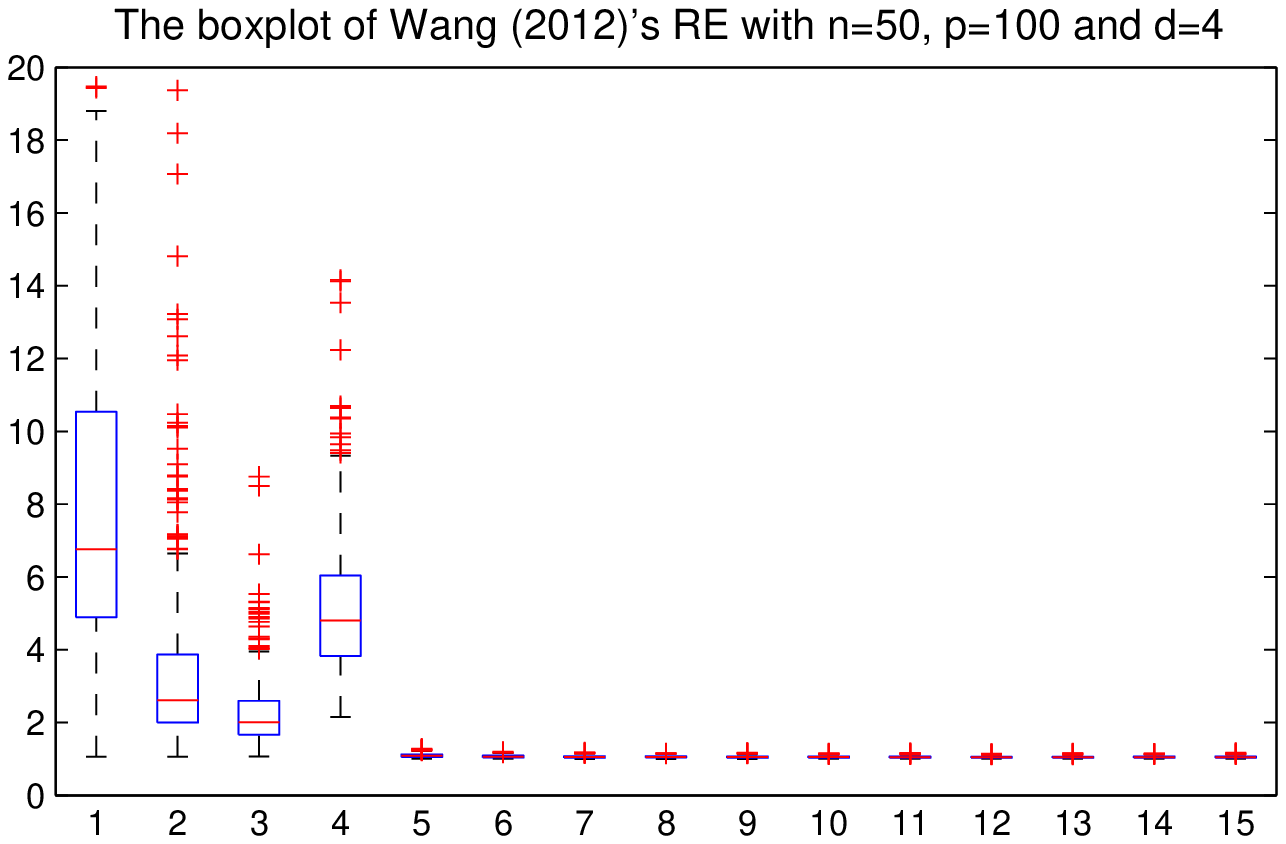}\\
  \caption{Boxplots of the first $15$ components for $\{\hat{s}^*_{j+1}+\frac{\log{(m)}}{5\sqrt{m}}\}/\{\hat{s}^*_j+\frac{\log{(m)}}{5\sqrt{m}}\}$ (TDRR),  $\{\hat{\lambda}_j+\frac{\log{n}}{10n}\}/\{\hat{\lambda}_{j+1}+\frac{\log{n}}{10n}\}$ (RRE),
 \{$\frac{\hat{\lambda}_j}{\hat{\lambda}_{j+1}}$\} (RE) and
  $\{\frac{n\sum_{l=1}^j\{\log(1+\hat{\lambda}_l)+\hat{\lambda}_l\}}
  {2\sum_{l=1}^p\{\log(1+\hat{\lambda}_l)+\hat{\lambda}_l\}}-\log(n)\frac{j(j+1)}{p}\}$ (BIC), for $j=1,\ldots,15$,  about the approximate factor model of Example~5.
  }\label{figure-factor}
\end{figure}

The findings from the results reported in Tables~\ref{factor-N} and~\ref{factor-T} are as follows. First,
TDRR can  accurately determine the number of common factors  at a high percentage (more than 90\%) no matter whether the components of $\mathbf{F}_t$ are correlated or not. Second, when the components of $\mathbf{F}_t$  are independent of  each other, the  three competitors also perform well and accurately estimate the number of common factors at  high percentages (more than 90\%). But when the components are highly  correlated, they seriously underestimate $d$, especially for either the small sample size $n$ or the small dimension $p$. This implies that our method is also robust against the correlations among the components of $\mathbf{F}_t$. Third, the comparison between  Tables~\ref{factor-N} and~\ref{factor-T} suggests that TDRR has a stable performance against  heavy tails. In contrast, the heavy tails of $\mathbf{F}_t$ have  negative impact for RRE, RE and BIC. To clearly demonstrate these  phenomena, we also present the boxplots in Figure~\ref{figure-factor}, which draw the results for the first fifteen components for different criteria when $n=50$ and $p=100$.

\begin{center}
Figure~\ref{figure-factor} about here
\end{center}
As the reasons why TDRR works well and the others do not can be very much similarly explained as showed by Figure~\ref{figure-reduction}, we then only give a very brief explanation here. Figure~\ref{figure-factor} shows that  the estimated numbers of common factors by RRE and RE are about $1$ and thus, the underestimation problem of these two methods is again confirmed. This is again because of the dominant effect of the largest eigenvalue over other nonzero eigenvalues.  Although BIC works better than  RRE and RE, it still underestimates the number $d$ of common factors when the covariance matrix is $\Sigma_2$.

Thus, for approximate factor models, TDRR is also an efficient method.

\subsection{Real data example: sufficient dimension reduction}

\begin{figure}
  \centering
  \includegraphics[width=7.250cm]{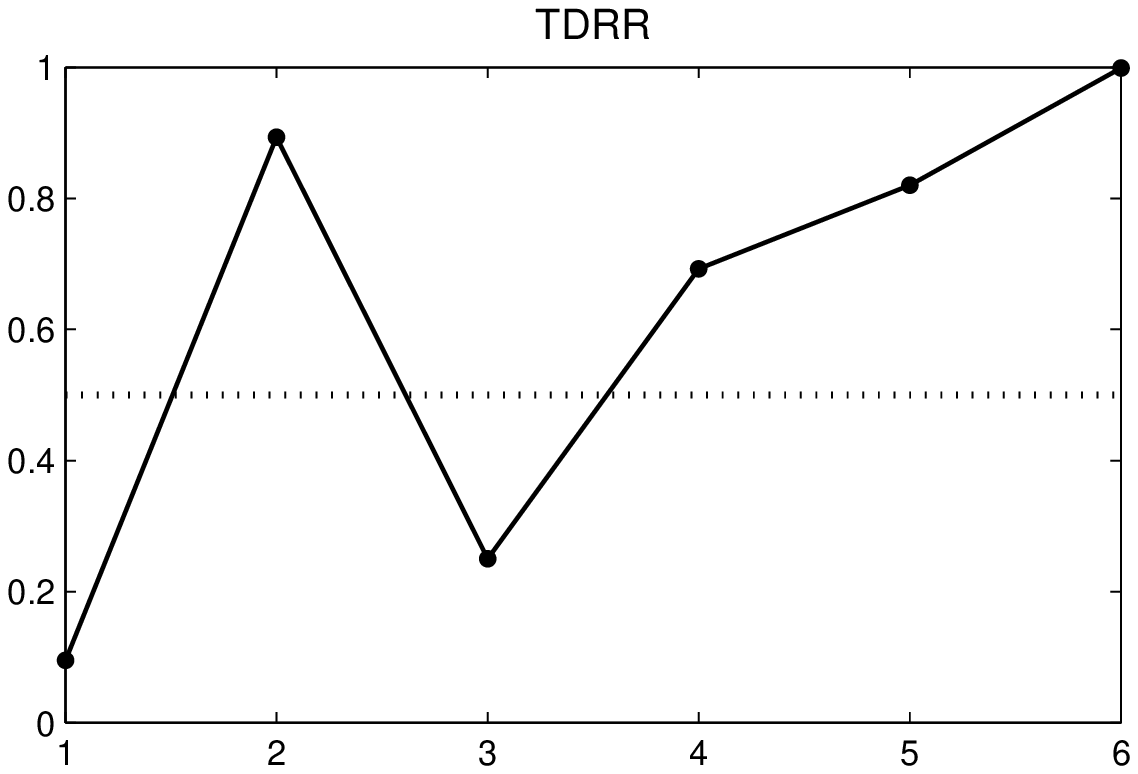}
  \includegraphics[width=7.250cm]{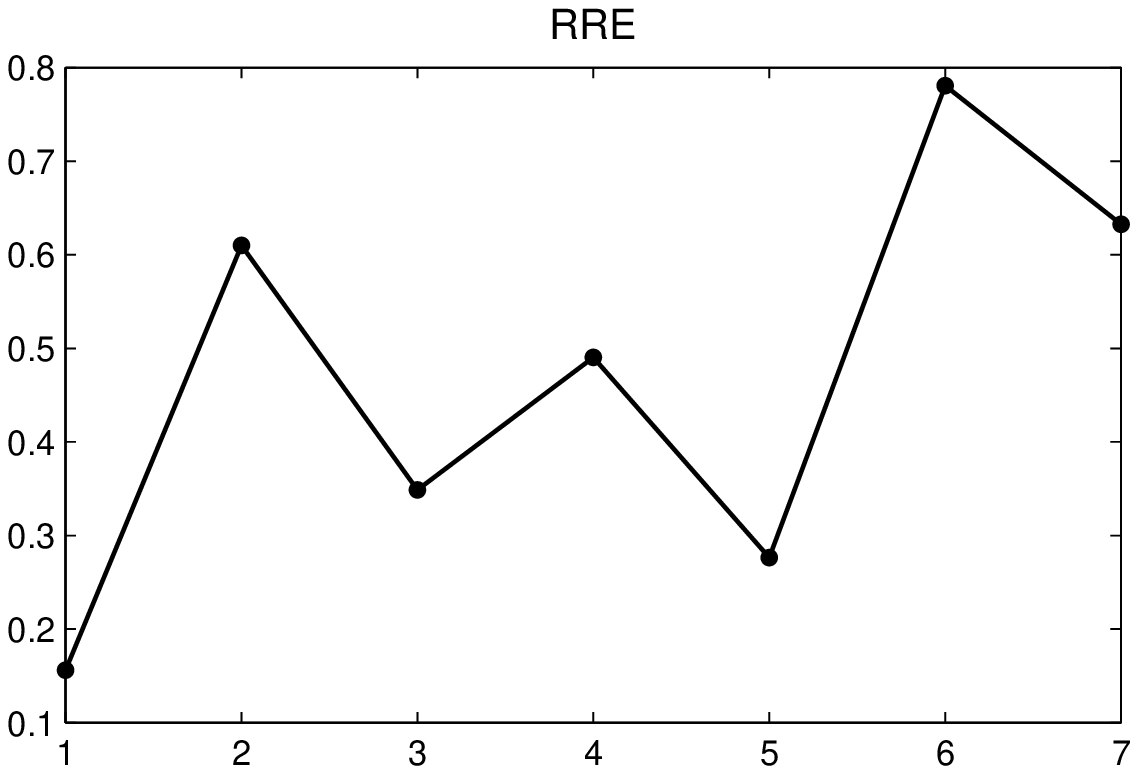}
  \includegraphics[width=7.250cm]{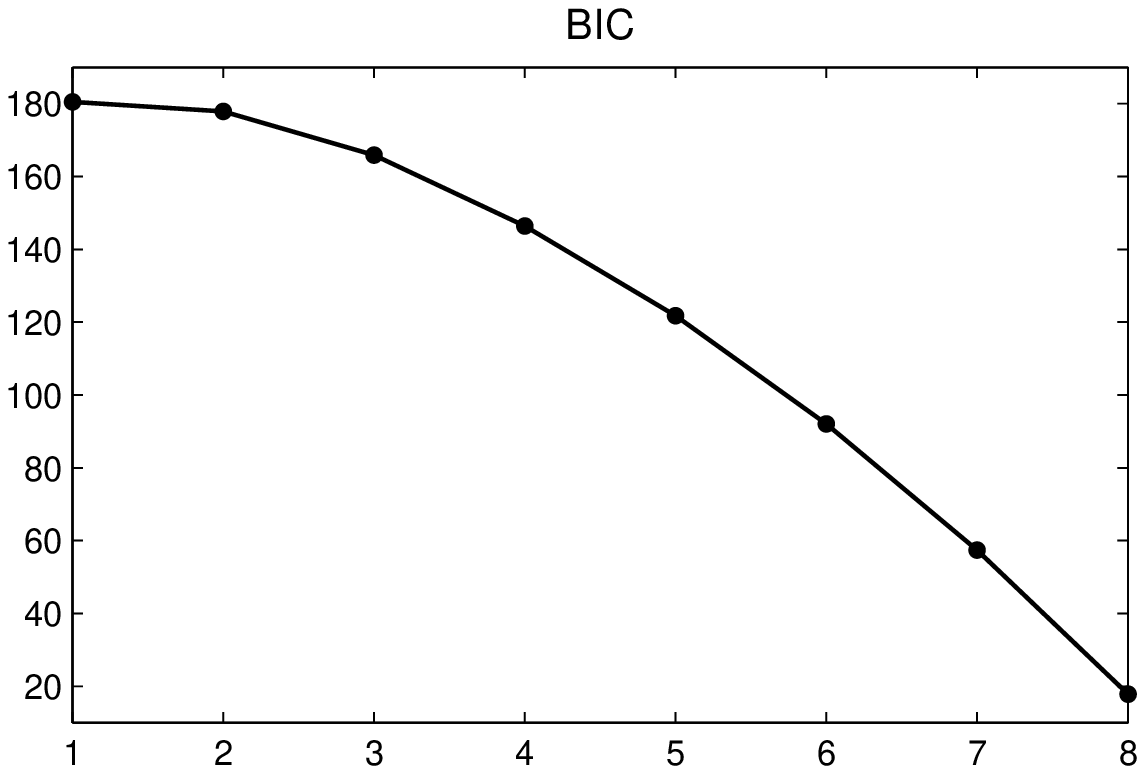}
  \includegraphics[width=7.250cm]{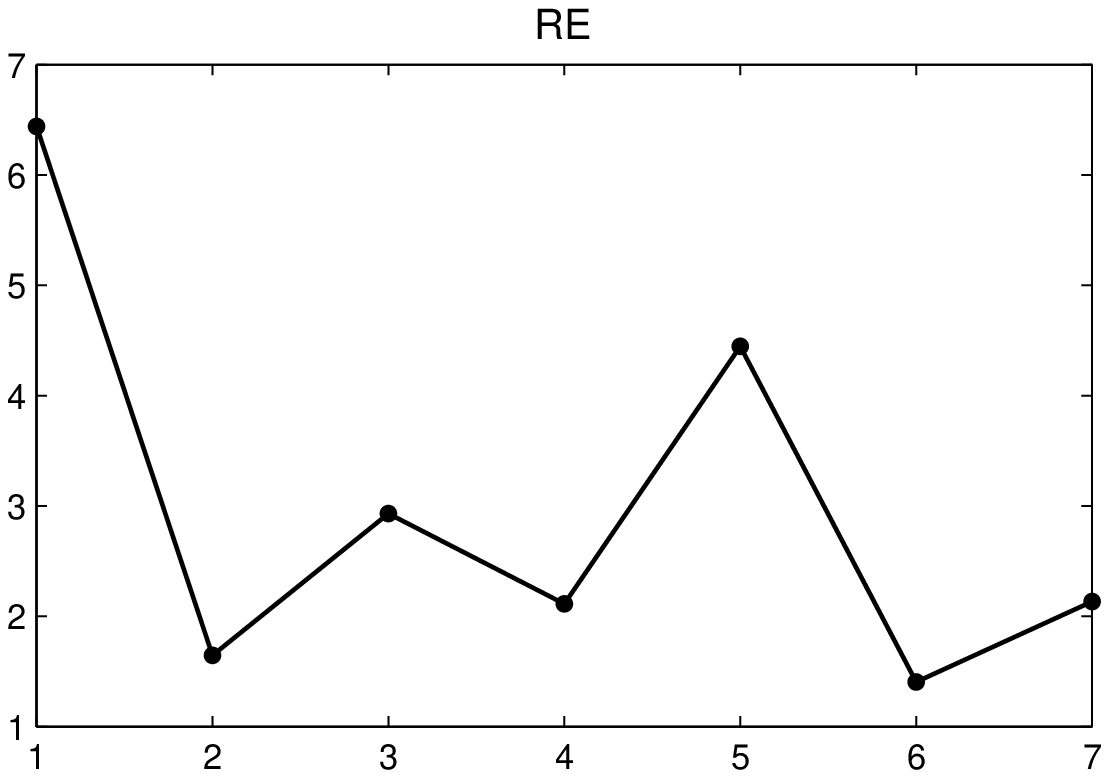}\\
  \caption{Plots of the components for $\{\hat{s}^*_{j+1}+\frac{\log{(m)}}{5\sqrt{m}}\}/\{\hat{s}^*_j+\frac{\log{(m)}}{5\sqrt{m}}\}$ (TDRR),  $\{\hat{\lambda}_j+\frac{\log{n}}{10n}\}/\{\hat{\lambda}_{j+1}+\frac{\log{n}}{10n}\}$ (RRE),
 \{$\frac{\hat{\lambda}_j}{\hat{\lambda}_{j+1}}$\} (RE) and
 $\frac{n\sum_{l=1}^j\{\log(1+\hat{\lambda_l})+\hat{\lambda}_l\}}
  {2\sum_{l=1}^p\{\log(1+\hat{\lambda_l})+\hat{\lambda}_l\}}-\sqrt{n}\frac{j(j+1)}{p}$ (BIC), for $j=1,\ldots,8$,  about the dimension reduction model of the Cars Data set.
  }\label{real-dimension}
\end{figure}

\begin{figure}
  \centering
  \includegraphics[width=15cm]{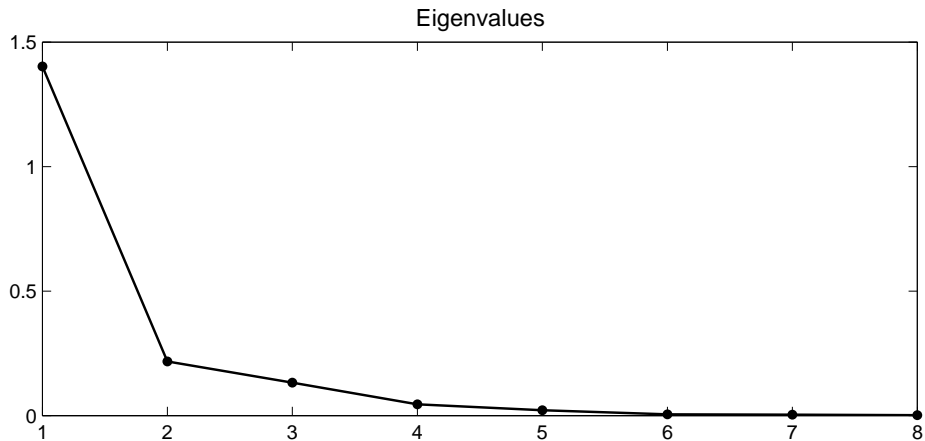}\\
  \caption{The plot of the estimated eigenvalues about the dimension reduction model of the Cars Data set.
  }\label{real-eigenvalues-dimension}
\end{figure}

Cars data set was used in the 1983 American Statistical Association Exposition, which is available at \url{http://archive.ics.uci.edu/ml/datasets/Auto+MPG}. This data set was analysed in Guo et al. (2015) and Xia (2007). There are 392 sample points and 8 variables, which
include the number of cylinders ($X_1$), engine displacement ($X_2$), horsepower ($X_3$), vehicle weight ($X_4$), time to accelerate from 0 to 60 mph ($X_5$), model year ($X_6$), origin of the car (1 = American, 2 = European, 3 = Japanese) and  miles per gallon ($Y$).
According to Xia (2007)'s suggestions,  we also take miles per gallon as the response variable and define two indicator variables because the origin of the car contains more than two categories. Let $X_7 = 1$ if a car is from America and 0 otherwise. Let $X_8 = 1$ if a car is from Europe and 0 otherwise. We standardize the predictors $X=(X_1,\cdots, X_8)^{\top}$. This is not a really high-dimensional problem. However, as we do not have prior information on the model structure,  the dimension $8$ is still regarded as high in nonparamatric estimation compared with a sample of size less than $400.$
Then we fit the data by using the following multi-index model with  an unknown number $q$ of indices:
\begin{eqnarray}\label{real}
 Y=G(B^{\top}X)+\epsilon.
\end{eqnarray}
As done before, we use the SIR-based DEE matrix $M_{n}$ and then determine the number $q$ by the four criteria we considered in the simulations.
\begin{center}
Figure~\ref{real-dimension} about here
\end{center}
\begin{center}
Figure~\ref{real-eigenvalues-dimension} about here
\end{center}

TDRR determines $\hat q=3$, while the other three competitors estimate $\hat q=1$. The results in Figure~\ref{real-dimension} suggest clearly that minimizers or maximizers of the criteria take the value of $1$. On the contrast, the thresholding strategy makes TDRR determine larger $\hat q$.
To interpret the reason why to happen these phenomena, we plot the estimated eigenvalues $\hat{\Lambda}=(1.4019,0.2177,0.1322,0.0451,0.0213,0.0048,0.0034,0.0016)$ in Figure~\ref{real-eigenvalues-dimension}. Clearly, the first eigenvalue  of value $1.4019 $ is very much larger than the second one of value $0.2177$ and the others, thus its dominating effect causes underdetermiantion. However, the fourth seems much smaller than the third, and thus, $\hat q=3$ seems more reasonable, which is determined by TDRR.  By the DIR-based DEE, we can obtain the three eigenvectors as the indices when TDRR is used and the corresponding first eigenvector as the index when the other competitors are applied. After that, we can estimate multi-index and single-index functions in the form of (\ref{real}) respectively to see what model fits data well:
\begin{eqnarray*}
\hat{G}(\hat{B}^{\top}_{\hat{q}}X) = \frac{\sum_{i=1}^nK_{\hat{q}h}(\hat{B}^{\top}_{\hat{q}}X-\hat{B}^{\top}_{\hat{q}}X_i)y_i}
{\sum_{i=1}^nK_h(\hat{B}^{\top}_{\hat{q}}X-\hat{B}^{\top}_{\hat{q}}X_i)}.
\end{eqnarray*}
where  $K_{\hat{q}h}=K(\cdot/h)/h^{\hat{q}}$ with $K(\cdot)$ being a $\hat{q}$-dimensional kernel function and $h$ being a bandwidth and $\hat{B}_{\hat{q}}$ is a sufficient dimension reduction estimate with an estimated structural dimension $\hat{q}$ of $q$.  We choose the product of $\hat{q}$ Quartic kernel function as $K(u) = 15/16(1 - u^2)^2$, if $|u| \leq 1$ and 0 otherwise and the bandwidth as $h= n^{-1/(4+\hat{q})}/4$. See Guo et al. (2015) for more details. We compute the residual sum of squares $RSS=\sum_{i=1}^n\{y_i-\hat{G}(\hat{B}^{\top}_{\hat{q}}x_i)\}^2/n$ To select the bandwidth $h$, we first use the cross validation. When $\hat q=1$, the $RSS$ value is $56.7057$ and when $\hat =3, 4$, the $RSS$' are about $37$. This suggests that the models with $\hat q=3$ and $4$ fit the data similarly.
However, the CV-based bandwidths seem not good choice as the corresponding RSS is still fairly large.  Then we manually select it by using grid points $l/20$ for $l=1, \cdots, 20$ within the interval (0.05, 1). The minimum value of  RSS is $47.6360$ when $\hat q=1$. When $\hat q=3$, the RSS value is quickly dropped down to $1.0170$. To check whether the model with $\hat q=3$ can well fit the data, we also check the RSS values with $\hat q=4$ and  all covariates. The values are respectively $0.8495$ and $1.7801$. These again indicate that $\hat q=3$ is a good choice. Although the RSS with $\hat q=4$ is slightly smaller, the model is less parsimonious in a nonparamatric setup and would make further analysis more difficult. When we use all covariates with$\hat q=8$, the RSS takes larger value showing that nonparametric estimation suffers from the curse of dimensionality.  Balancing between model parsimoniousness and model fitting, $\hat q=3$ is a good choice. 

In conclusion, the new method can be immune to the dominating effect by the largest eigenvalues in dimensionality determination whereas the other competitors  can not.

\subsection{Real data example: approximate factor model}
In this subsection, we apply our method to analyze the dataset, which was used to determine the number of common factors in the approximate factor models in Stock and Watson (2005) and Caner and Han (2014).  The U.S. macroeconomic data in this dataset, spanning the period of 1960.1-2003.12. The covariates are transformed to achieve stationarity  and then outliers are adjusted as described in Appendix A of Stock and Watson (2005). It then contains 132 covariates and 526 observations. The more details about the covariates can be referred to Stock and Watson (2005).

\begin{figure}
  \centering
  \includegraphics[width=15cm]{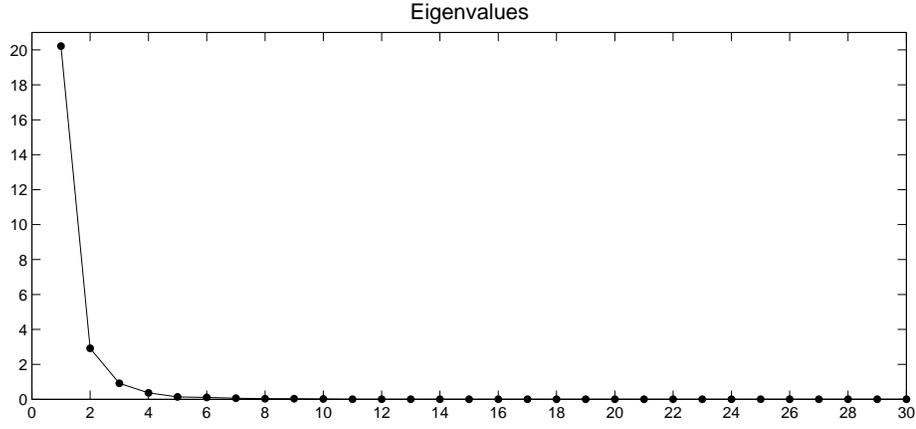}\\
  \caption{The plot of the first $30$ estimated eigenvalues about the approximate factor model of the real data.
  }\label{real-eigenvalues}
\end{figure}

Similarly as Figure~\ref{real-eigenvalues-dimension}, we plot the first 30 estimated eigenvalues in Figure~\ref{real-eigenvalues}.
The first eigenvalue of  $20.2168$ is much larger than the second one of  $2.9217$,  thus the  dominating effect also exists.

\begin{center}
Figure~\ref{real-eigenvalues} about here
\end{center}
The first $30$ values in the criteria are plotted in Figure~\ref{real-factor} below.

\begin{center}
Figure~\ref{real-factor} about here
\end{center}
The plots obviously indicate that by TDRR, there are three local minima at $d=4, 7$ and $10$ and the resultant estimator $\hat d=10$, BIC estimates  $\hat d=4$, and
both  RRE  and RE  only determine the number $d$ to be $1$.  When  an optimization method suggested in Bai and Ng (2002) was adopted by Stock and Watson (2005), which  is also a kind of information criterion requiring heavy computational workload, the number was determined to be $9$ that is close to the number TDRR determines for a conservative consideration. A relevant reference is Caner and Han (2014). Therefore, underdetermination to make  working models too parsimonious by the three competitors is an indispensable difficulty.

\section{Concluding remarks}
In this paper, we propose a novel approach to determine the structural dimension of the central (mean) subspaces in sufficient dimension reduction. It is a generic method and can also be applied to determine the number of common factors in approximate factor models even when the dimension is ultra-high. The new method can largely solve the longstanding problem of underdetermination for all existing eigendecomposition-based methods in the literature. The numerical studies show its usefulness and superiority to the competitors.
This method could also be useful to determine the number of principal components for principal component analysis in the high dimensional analysis and other types of data such as functional data (e.g. Li and Hsing, 2010), longitudinal data (e.g. Jiang, et al., 2014; Bi and Qu, 2015) and tensor data. The relevant researches are ongoing.

Another issue about overdetermination in some scenarios as we have seen in the numerical studies. To obtain a more parsimonious working model without losing important covariates or linear combinations of covariates, we may consider a further selection  after our method is implemented. A naive idea would be that we determine the dimensionality within the selected linear combinations or covariates that are regarded as all new covariates. Note that the number $\tilde q$ of the new covariates should be much smaller than the original dimension $p$, the second round of selection would be efficient. However, based on  our initial numerical study, it seems not to work well as the second round of selection would not very efficiently reduce $\tilde q$ to $q$ if $\tilde q> q$. On the other hand, if $\tilde q$ is already smaller than $q$, the new round of selection would not have a further reduction either. In other words, a second round of selection seems to only slightly change the selection by our method. Thus, how to significantly improve our method in this scenario deserves a further study.

\section{Appendix: Proofs}
\subsection{Regularity Conditions}
To prove  Theorem~\ref{theorem-seq} about the local models, we need the following conditions.
\begin{itemize}
\item [A3.] The linearity condition on the distribution of $X$ is satisfied,  namely, the conditional mean $E(X|C^{\top}X)$
 is linear in $C^{\top}X$,  where $C \in \mathbb{R}^{p\times q}$ is any basis matrix for the subspace $\rm{Span}(B_1, B)$.

\item [A4.]  $\tilde{M}$ has $q$ nonzero eigenvalues, where $\tilde{M}(t)=\left\{\tilde{m}(t)m(t)^{\top}+m(t)\tilde{m}(t)^{\top}\right\}$ with $\tilde m(t)=E\big (-g(B^{\top}X)f_{Y|X}(t)X\big )$ and $m(t)=E\{(X-E(X))I(Y\leq t)\}$ and  $f_{Y|X}(\cdot)$ stand for  the conditional density function of $Y$ given $X$. 
\end{itemize}

To prove  Theorem~\ref{theorem-factor} for the approximate factor model, the same conditions in Wang (2012) are designed as follows:
\begin{itemize}
\item [B1.] The common factor matrix $\mathbf{F}$ and the idiosyncratic factor matrix $\mathbf{U}$ are normally distributed. Additionally,  there exists some positive constant $\sigma_{min}$ satisfying  $\min_{1\geq j \geq p}\sigma^2_{j} \sigma^2_{min} >0$, where $\sigma^2_{j}={\rm{var}} (u_{ij})$.
\item [B2.] Both $B_j$ and $\sigma^2_{j}$ admit $p^{-1}\mathbf{B}^{\top}\mathbf{B}=p^{-1}\sum_{j=1}^pB_jB^{\top}_j=\Sigma_{B}+O_p(p^{-1/2})$ and $p^{-1}\sum_{j=1}^p\sigma^2_{j}=\sigma^2_{0}+O_p(p^{-1/2})$ with $\Sigma_{B}$ being a $d\times d$ positive definite matrix and $\sigma^2_{0}$ being a positive constant.
\item [B3.] The number $d$ of common factors is fixed and $\log(p)=O_p(n^h)$ with $0<h<1$.
\end{itemize}

\subsection{Proofs}
\begin{lemma}\label{lemma1}
If $A$ and $B$ are $k\times k$ symmetric matrixes, then
\begin{eqnarray}
\lambda_{i+j-1}(A+B)\leq \lambda_{i}(B)+\lambda_{j}(B)~{\rm{for}}~i+j-1\leq k,
\end{eqnarray}
where $\lambda_{i}(C)$ denote the $i$-th largest eigenvalue of the matrix $C$.
 \end{lemma}
\textbf{Proof of Lemma~\ref{lemma1}.} See Ahn and Horenstein (2013).

\begin{lemma}\label{lemma2}
If $A$ and $B$ are $k\times k$ positive semi-definite matrixes, then
\begin{eqnarray}
\lambda_{i}(A)\leq \lambda_{i}(A+B)~{\rm{for}}~i \leq k,
\end{eqnarray}
where $\lambda_{i}(C)$ denote the $i$-th largest eigenvalue of the matrix $C$.
\end{lemma}
\textbf{Proof of Lemma~\ref{lemma2}.} See Ahn and Horenstein (2013).

\textbf{Proof of Theorem~\ref{theorem0}. } 
Under Condition A1, we have $M_{n}-M=O_p(n^{-1/2})$.
Following the similar arguments used in Zhu and Ng (1995) or Zhu and Fang (1996), we can prove   the root-$n$ consistency of the eigenvalues of $M_{n}$, namely, $\hat{\lambda}_{i}- \lambda_{i} =O_p(n^{-1/2})$. Recall the definition that
\begin{eqnarray*}
s_{j}=\frac{\lambda_{j}}{1+\lambda_{j}}\ {\rm{\ for}}\  1\leq j < p,
\end{eqnarray*}
where $ \lambda_{1}\geq \ldots \geq\lambda_{p}\ge 0$ denote the eigenvalues of the target matrix $M$.
Also recall the definition of  ${s}^*_j$ for all $1\le j<p-1$:
\begin{eqnarray}\label{lambdastar}
s^*_j=
\left\{\begin{array}{ll}
\left(\frac{s_{j}}{s_{j+1}}\right)^2-1\ge 0,
 & {\rm{\ for}}\  1\leq j < q,\\
+\infty, & {\rm{\ for}}\   j = q,\\
0, & {\rm{\ for}}\  q+1 \leq j \leq p-1.
\end{array}\right.
\end{eqnarray}
We  now prove that these are the limits of $\hat{s}^*_{j} $ for all $1\le j<p$.
Because the function $x/(1+x)$ is continuous and strictly monotone about $x\ge 0$, we have $\hat s_{j}-s_{j}=O_p(1/\sqrt{n})$. Notice that $s_j^2=0$ for $q+1\leq j \leq p$, then we have
\begin{eqnarray*}
\hat s^2_{j} =\left\{\begin{array}{ll}
s^{2}_{j}+O_p(1/\sqrt{n}), & {\rm{\ for}}\  1\leq j \le q,\\
O_p(1/n), & {\rm{\ for}}\  q+1\leq j \leq p.
\end{array}\right.
\end{eqnarray*}
Recall the definition in  (\ref{tilderatio}) that  $\hat{s}^*_q=\frac{\hat s^2_{q}+c_{1n}}{\hat s^2_{q+1}+c_{1n}}-1=O_p(1/c_{1n})$. Then it is easy to justify that as $n\to \infty$, $\hat{s}^*_{q}  \rightarrow  +\infty$.
For $j$ with $ 1\leq j < q$,  we have
\begin{eqnarray*}
\frac{\hat s^2_{j}+c_{1n}}{\hat s^2_{j+1}+c_{1n}}-1=\frac{s^2_{j}+O_p(1/\sqrt{n})+c_{1n}}{s^2_{j+1}+O_p(1/\sqrt{n})+c_{1n}}-1=
\frac{s^2_{j}}{s^2_{j+1}}-1+O_p(\max\{1/\sqrt{n},c_{1n}\}).
\end{eqnarray*}
When $ q+1 \leq j \leq p-1$,  we have
\begin{eqnarray*}
\frac{\hat s^2_{j}+c_{1n}}{\hat s^2_{j+1}+c_{1n}}-1=\frac{\hat s^2_{j}-\hat s^2_{j+1} }{\hat s^2_{j+1}+c_{1n}}=
\frac{O_p(1/n)}{O_p(1/n)+c_{1n}}=O_p(\frac{1}{nc_{1n}}).
\end{eqnarray*}
Altogether, we derive
\begin{eqnarray}\label{lambdastar1}
\hat{s}^*_{j} =\left\{\begin{array}{ll}
 s^*_{j}+O_p(\max\{\frac1{\sqrt{n}},c_{1n}\}), & \rm{\ for}\  1\leq j < q,\\
O_p(\frac{1}{c_{1n} }), & j=q,\\
O_p(\frac{1}{nc_{1n} }), &\rm{\ for}\  q+1 \leq j \leq p-1.
\end{array}\right.
\end{eqnarray}
Therefore, $s^*_{j}$ in (\ref{lambdastar}) are the limits of $\hat{s}^*_{j}$. The ratios in the second round satisfy that as $n\to \infty$,
for any $1\le j<q-1$
\begin{eqnarray*}
\frac{\hat{s}^*_{j+1}+c_{2n}}{\hat{s}^*_{j}+c_{2n}}
\rightarrow \frac{{s}^*_{j+1}}{{s}^*_{j}},
\end{eqnarray*}
for $j=q-1$
\begin{eqnarray*}
\frac{\hat{s}^*_{q}+c_{2n}}{\hat{s}^*_{q-1}+c_{2n}}
\rightarrow \frac{+\infty}{{s}^*_{q-1}}=+\infty,
\end{eqnarray*}
for  $j=q$
\begin{eqnarray*}
\frac{\hat{s}^*_{q+1}+c_{2n}}{\hat{s}^*_{q}+c_{2n}}
=\frac{c_{2n}+O_p(\frac{1}{nc_{1n} })}{\hat{s}^*_{q}+c_{2n}}
\rightarrow \frac{0}{+\infty} = 0,
\end{eqnarray*}
and  for any $j>q$, $s^*_{j}=0$,
\begin{eqnarray*}
\frac{\hat{s}^*_{j+1}+c_{2n}}{\hat{s}^*_{j}+c_{2n}}
=\frac{c_{2n}+O_p(\frac{1}{nc_{1n} })}{c_{2n}+O_p(\frac{1}{nc_{1n} })} \rightarrow  1
\end{eqnarray*}
by the condition $c_{1n}c_{2n}n \rightarrow  \infty$.
Altogether, we have
\begin{eqnarray}\label{ratio1}
\lim_{n\rightarrow \infty}\frac{\hat{s}^*_{j+1}+c_{2n}}{\hat{s}^*_{j}+c_{2n}}=\left\{\begin{array}{ll}
\frac{{s}^*_{j+1}}{{s}^*_{j}} >0, & {\rm{\ for}}\   j <q-1,\\
+\infty > \tau, & {\rm{\ for}}\   j = q-1,\\
0 < \tau, & {\rm{\ for}}\   j = q,\\
1>\tau, & {\rm{\ for}}\  q+1 \leq j \leq p-2.
\end{array}\right.
\end{eqnarray}
Therefore, we can conclude that as $n\rightarrow \infty$, $\hat{q}= q$ with a  probability going to $1$. \hfill$\Box$

\textbf{Proof of Theorem~\ref{theorem2.2}.}
Under the condition in Theorem~\ref{theorem2.2}, we have $||M_n-M||=O_p(p^{1/2}n^{-1/2})$.
By the similar arguments in Zhu and Ng (1995) or Zhu and Fang (1996), it is proved that  $\hat{\lambda}_{i}- \lambda_{i} =O_p(p^{1/2}n^{-1/2})$. Here we use the same notations in the justification in Theorem~\ref{theorem0}.
The continuous function $x/(1+x)$ results in  $\hat s_{j}-s_{j}=O_p(p^{1/2}n^{-1/2})$. This leads to
\begin{eqnarray*}
\hat s^2_{j} =\left\{\begin{array}{ll}
s^{2}_{j}+O_p(p^{1/2}n^{-1/2}), & {\rm{\ for}} \  1\leq j \leq q,\\
O_p(p/n), &{\rm{\ for}} \  q+1\leq j \leq p-1.
\end{array}\right.
\end{eqnarray*}
When taking $c_{1n}c_{2n}n/p \rightarrow  \infty$, the limits of the ratios $\lim_{n\rightarrow \infty}\frac{\hat{s}^*_{j+1}+c_{2n}}{\hat{s}^*_{j}+c_{2n}}$ have the same properties as those of (\ref{ratio1}).  Thus, as $n\rightarrow \infty$,  $\hat{q}= q$  with a probability going to one.
 \hfill$\Box$

\textbf{Proof of Theorem~\ref{theorem-seq}. }
By the SIR-based DEE procedure, we have
\begin{eqnarray*}
M(t)=\Sigma^{-1}{\rm{var}}\left[E\{X|Z(t)\}\right]=\Sigma^{-1}(\nu_{t1}-\nu_{t0})(\nu_{t1}-\nu_{t0})^{\top}
p_t(1-p_t),
\end{eqnarray*}
where $\Sigma$ is the covariance matrix of $X$, $\nu_{tj}=E\{X|Z(t)=j\}$ for $j=0$ and 1, and $p_t=E\{I(Y\leq t)\}$. Further note that:
\begin{eqnarray*}
\nu_{t1}-\nu_{t0}&=&\frac{E\{XI(Y\leq t)\}}{p_t}-\frac{E\{XI(Y>t)\}}{1-p_t}\\
&=&\frac{E\{XI(Y\leq t)\}-E(X)E\{I(Y\leq
t)\}}{p_t(1-p_t)}.
\end{eqnarray*}
Applying Lemma~1 in Guo et al. (2015),  $M(t)$ can also be rewritten as:
\begin{eqnarray*}
M(t)&=&\Sigma^{-1}\Big(E[\{X-E(X)\}I(Y\leq t)]\Big)\Big(E[\{X-E(X)\}I(Y\leq t)]\Big)^{\top}\\
&=:&\Sigma^{-1} m(t)m(t)^{\top}=\Sigma^{-1}L(t),
\end{eqnarray*}
where $m(t)=E\{(X-E(X))I(Y\leq t)\}$. Then $m(t)$ and $L(t)$ can be respectively estimated by
\begin{eqnarray}\label{LN}
m_n(t)&=&n^{-1}\sum_{i=1}^n(X_i-\bar{X})I(y_i\leq t), \nonumber\\
L_{n}({t})&=&m_{n}(t)m^{\top}_{n}(t),
\end{eqnarray}
and then  $M(t)$ can be estimated by $M_{n}(t) =\hat{\Sigma}^{-1}L_{n}(t). $
Since the responses under the local models in (\ref{sequence}) are related to the sample size $n$, we rewrite the responses under $C_n=0$ and $C_n \neq 0$ as $Y$ and $Y_n$ respectively and rewrite $y_{in}$ as the responses at the sample level. Here $Y_n$ is the true response and $Y$ can be regarded as the limit to whom $Y_n$ converges.
Thus, from Zhu et al. (2010), we have
\begin{eqnarray*}
&&\frac{1}{n}\sum_{i=1}^nx_iI(y_{in}\leq t)-E\{XI(Y\leq t)\}\\
&=&\frac{1}{n}\sum_{i=1}^n[x_iI(y_{in}\leq t)-E\{xI(Y_n\leq t)\}]+E\{XI(Y_n\leq t)\}-E\{XI(Y\leq t)\}\\
&=&O_p(n^{-1/2})+E\{XI(Y_{n} \leq t)\}-E\{XI(Y \leq t)\}.
\end{eqnarray*}
Further,
\begin{eqnarray*}
E\{XI(Y_n\leq t)\}-E\{XI(Y\leq t)\}=E\Big[X\{P(Y_n\leq t|X)-P(Y\leq
t|X)\}\Big].
\end{eqnarray*}
Under the local models in (\ref{sequence}), because  $Y_n = Y+C_n g(B^{\top}X)$, we have for all $t$,
\begin{eqnarray}\label{Taylor}
P(Y_{n} \leq t|X)-P(Y \leq t|X)
&=&F_{Y|X}\{t-C_ng(B^{\top}X)\}-F_{Y|X}(t)\nonumber\\
&=&-C_ng(B^{\top}X)f_{Y|X}(t)+O_p(C_n^2),
\end{eqnarray}
where $f_{Y|X}(\cdot)$ and $F_{Y|X}(\cdot)$ respectively stand for  the conditional density and  distribution function of $Y$ given $X$.

\textit{ Proof of Part (I).}
The formula (\ref{Taylor}) implies that
$$n^{-1}\sum_{i=1}^nx_iI(y_{ni} \leq t)-E\{XI(Y\leq
t)\}=O_p (\max\{C_n, n^{-1/2} \}).$$
It is easy to derive that
\begin{eqnarray*}
m_n(t)-m(t)&=& O_p (\max\{C_n, n^{-1/2} \}),\\
L_n(t)-L(t)&=&O_p (\max\{C_n, n^{-1/2} \}).
\end{eqnarray*}
These results yield that
\begin{eqnarray*}
M_n(t)-M(t)=O_p (\max\{C_n, n^{-1/2} \}).
\end{eqnarray*}
Finally,  similar to the argument used for proving Theorem 3.2 of Li et al. (2008),
$M_{n}-M=O_p (\max\{C_n, 1/\sqrt{n} \})$. When $C_n =O(1/\sqrt{n}) $, $O_p (\max\{C_n, 1/\sqrt{n} \})=O_p(1/\sqrt{n})$. Thus, the eigenvalues of $M_n$ converge to the corresponding eigenvalues of $M$ at the rate of $O_p(1/\sqrt{n}).$ Applying the same argument used for proving  Theorem~\ref{theorem0}, we  conclude that as $n\rightarrow \infty$, $\hat{q}= q_1$ with a probability going to one.

\textit{ Proof of Part (II).}
Similarly as (\ref{Taylor}), we consider the second derivative of $F_{Y|X}\{t-C_ng(B^{\top}X)\}-F_{Y|X}(t)$ about $t$ to obtain that
\begin{eqnarray*}
&&F_{Y|X}\{t-C_ng(B^{\top}X)\}-F_{Y|X}(t)\\
&=&-C_ng(B^{\top}X)f_{Y|X}(t)+(C_ng(B^{\top}X))^2f'_{Y|X}(t)/2+O_p(C_n^2).
\end{eqnarray*}
Let $\tilde m(t)=E\big \{-g(B^{\top}X)f_{Y|X}(t)X\big \}$ and $m(t)=E\big [\{X-E(X)\}I(Y\leq t)\big ]$. By the definition of $L_n$ in (\ref{LN}), we have
\begin{eqnarray*}
M_{n}(t)&=&  \hat{\Sigma}^{-1}L_{n}(t)\\
&=&M(t)+C_n\Sigma^{-1}\left\{\tilde{m}(t)m(t)^{\top}+m(t)\tilde{m}(t)^{\top}\right\}
+O_p(\max\{C^2_n,n^{-1/2}\})\\
&=:&M(t)+C_n\tilde{M}(t)+O_p(\max\{C^2_n,n^{-1/2}\})
\end{eqnarray*}
where $\tilde{M}(t)=\left\{\tilde{m}(t)m(t)^{\top}+m(t)\tilde{m}(t)^{\top}\right\}$.

When $C_n=0$, as $Y \hDash X|B^{\top}_1X$,  $f_{Y|X}(t)$ and $f'_{Y|X}(t)$ are the functions of $(B^{\top}_1X, t)$, we then define them as $f_{Y|B^{\top}_1X}( t)$ and $f'_{Y|B^{\top}_1X}( t)$. Under the linearity condition A3, $E(X|B^{\top}_1X,B^{\top}X)=P_{(B_1,B)}X$ where  $P_{(B_1,B)}$ denotes the projection operator of the subspace $\rm{Span}(B_1, B)$. 
Then
\begin{eqnarray*}
\tilde m(t)&=& E\{-g(B^{\top}X)f_{Y|X}(t)X\}=E\{-g(B^{\top}X)f_{Y|B^{\top}_1X}(t)X\}\\
&=&E\{-g(B^{\top}X)f_{Y|B^{\top}_1X}(t)E(X|B^{\top}_1X,B^{\top}X)\}\\
&=&P_{(B_1,B)}E\{-g(B^{\top}X)f_{Y|B^{\top}_1X}(t)X\}\in \rm{Span}(B_1, B),
\end{eqnarray*}
Similarly, we can prove that $m(t) \in \rm{Span}(B_1)$. Thus, $\tilde{M}(t)\in \rm{Span}(\tilde{M})$.


Following the argument for proving Theorem 3.2 in Li et al. (2008), we can derive that  $M_n(t)=M(t)+ C_n \tilde{M}(t)+O_p(\max\{C^2_n,n^{-1/2}\})$ uniformly over all $t$. Thus,
$M_n=M+ C_n \tilde{M}+O_p(\max\{C^2_n,n^{-1/2}\})$, where $\tilde{M}=E\{\tilde{M}(\tilde{Y})\}$ with $\tilde{Y}$ being  an independent copy of $Y$.

Under  the local models in (\ref{sequence}),  Condition A4 yields that $\rm{Span}(M)=\rm{Span}(B_1)$ and  $\rm{Span}(\tilde{M})=\rm{Span}(B_1, B)$ and then $\rm{Span}(M) \subseteq \rm{Span}(\tilde{M})$.  Then $\rm{rank}(\tilde M)=\rm{rank}(M+ C_n \tilde{M})=q$. The number of  nonzero eigenvalues of $M+ C_n \tilde{M}$ is the same as the dimension of the space $\rm{Span}(\tilde{M})$.
According to the definition of the local  models in (\ref{sequence}), let $\lambda_{11}\ge \ldots  \ge\lambda_{1q_1}>0=\lambda_{1(q_1+1)}=\ldots= \lambda_{1p}=0$ and $\lambda_{21}\ge\ldots \ge \lambda_{2q}>0=\lambda_{2(q+1)}=\ldots= \lambda_{2p}=0$ respectively denote the eigenvalues of the matrixes $M$ and $\tilde{M}$.  Further,
 as the eigenvalues of $M+ C_n \tilde{M}$ depends on the sample size $n$,
 define $\lambda_{1,n}\ge \ldots \ge \lambda_{q,n}>0=\lambda_{(q+1,n)}=\ldots= \lambda_{p,n}=0$  the eigenvalues of $M+ C_n \tilde{M}$.

Using Lemma~\ref{lemma1} and Lemma~\ref{lemma2}, we have
\begin{eqnarray*}
\left\{
\begin{array}{ll}
\lambda_{1i}\leq \lambda_{i,n}\leq \lambda_{1i}+C_n\lambda_{21}, & {\rm{for}} \ i=1,\ldots, q_1; \\
C_n\lambda_{2(q_1+1)}\leq \lambda_{i,n}\leq C_n\lambda_{21}, & {\rm{for}} \ i=q_1+1,\ldots, q;\\
\lambda_{i,n}=0,& {\rm{for}} \ i=q+1,\ldots, p.\\
\end{array}\right.
\end{eqnarray*}
Applying the mean value theorem, there exist some constants ${\mu_i}$ for $i=1, \ldots, q$ satisfying that:
\begin{eqnarray}\label{mu1}
\lambda_{i,n}=\left\{
\begin{array}{ll}
\lambda_{1i}+C_n\mu_{i}, & {\rm{for}} \ i=1,\ldots, q_1, \\
C_n\mu_{i}; & {\rm{for}} \ i=q_1+1 ,\ldots, q,\\
0,& {\rm{for}} \ i=q+1,\ldots, p,\\
\end{array}\right.
\end{eqnarray}
 where  $\mu_{1}\ge \ldots \ge \mu_{q_1}\ge 0$ are  nonnegative constants and $\mu_{q_1+1}\ge \ldots \ge \mu_{q}>0$ are some positive constants.

Similarly as that in  Zhu and Ng (1995) and Zhu and Fang (1996), we can prove that $\hat{\lambda}_{i}- \lambda_{i,n} =O_p(\max\{C^2_n,n^{-1/2}\})$.

Now we turn to compute the TDRR criterion.
From the definition of $\hat{s}^*_j$ in (\ref{tilderatio}), we want to show the following are their limits:
\begin{eqnarray}\label{star1}
s^*_j=
\left\{\begin{array}{ll}
\left(\frac{s_{j}}{s_{j+1}}\right)^2-1 >0,
 & {\rm{\ for}}\  1\leq j \le q_1-1,\\
+\infty, & {\rm{\ for}}\   j = q_1,\\
\tilde{\mu}^2_{j}-1 \ge 0,
 & {\rm{\ for}}\ q_1+1 \leq j \le q-1,\\
 +\infty, & {\rm{\ for}}\   j = q,\\
0, &{\rm{\ for}}\  q+1 \leq j \leq p-1,
\end{array}\right.
\end{eqnarray}

where $s_{j}=\frac{\lambda_{1j}}{1+\lambda_{1j}},\ \rm{\ for}\  1\leq j \le q_1,$ and $\tilde{\mu}_{j}=\mu_j/\mu_{j+1}\ge 1$ with $\mu_j$ defined in (\ref{mu1}).

For $q_1+1 \leq j \le q-1$, we have
$$\hat{s}_{j} = \frac{\hat{\lambda}_{j}}{\hat{\lambda}_{j}+1}=\frac{C_n\mu_j+O_p(\max\{C^2_n,n^{-1/2}\})}{1+C_n\mu_j+O_p(\max\{C^2_n,n^{-1/2}\})}=C_n\mu_j+o_p(C_n).$$
As  $c_{1n}=o_p(C^2_n)$, we have
\begin{eqnarray*}
\hat{s}^*_j+1&=&\frac{\hat{s}^2_{j}+c_{1n}}{\hat{s}^2_{j+1}+c_{1n}}
=\frac{(C_n\mu_j+o_p(C_n))^2+c_{1n}}{(C_n\mu_{j+1}+o_p(C_n))^2+c_{1n}}\\
&=&\frac{C^2_n\mu^2_j+o_p(C^2_n)+c_{1n}}{C^2_n\mu^2_{j+1}+o_p(C^2_n)+c_{1n}}
=\frac{\mu^2_{j}}{\mu^2_{j+1}}+o_p(1).
\end{eqnarray*}
Thus, for $q_1+1 \leq j < q-1$, $\hat{s}^*_j=\frac{\mu^2_{j}}{\mu^2_{j+1}}-1+o_p(1)$.
Consider the term $\hat{s}^*_q$. We have
\begin{eqnarray*}
\hat{s}^*_q+1&=&\frac{\hat{s}^2_{q}+c_{1n}}{\hat{s}^2_{q+1}+c_{1n}}
=\frac{C^2_n\mu^2_j+o_p(C^2_{1n})+c_{1n}}{O_p(\max\{C^4_n,n^{-1}\})+c_{1n}}
\rightarrow \infty,
\end{eqnarray*}
then  $\hat{s}^*_q\rightarrow \infty$. Altogether, $\hat{s}^*_j\rightarrow s^*_j$ defined in (\ref{star1}).

We now prove that, in probability,
\begin{eqnarray*}
\lim_{n\rightarrow \infty}\frac{\hat{s}^*_{(j+1)}+c_{2n}}{\hat{s}^*_{j}+c_{2n}}=\left\{\begin{array}{ll}
C_j>0, & \rm{\ for}\   1\le j \le q_1-1,\\
0, & \rm{\ for}\   j = q_1,\\
C_j>0, & \rm{\ for}\   q_1+1\le j \le q-2,\\
+\infty, & \rm{\ for}\   j = q-1,\\
0 < \tau, & \rm{\ for}\   j = q,\\
1>\tau, &\rm{\ for}\  q+1 \leq j \leq p-2.
\end{array}\right.
\end{eqnarray*}

As the proving argument for the ratios with different $j$ is similar, we only present the details for  the $j=q$  and $j>q$ case. Since $C^4_n/c_{1n}\rightarrow 0$, one derives
\begin{eqnarray*}
\frac{\hat{s}^*_{(q+1)}+c_{2n}}{\hat{s}^*_{q}+c_{2n}}
\rightarrow \frac{0}{\infty} =0<\tau.
\end{eqnarray*}
For any $j>q$,  we have
$$\hat{s}_{j} = \frac{\hat{\lambda}_{j}}{\hat{\lambda}_{j}+1}=\frac{O_p(\max\{C^2_n,n^{-1/2}\})}{1+O_p(\max\{C^2_n,n^{-1/2}\})}=O_p(\max\{C^2_n,n^{-1/2}\}).$$
Therefore, it is easy to see that
\begin{eqnarray*}
\hat{s}^*_j=\frac{\hat{s}^2_{j}+c_{1n}}{\hat{s}^2_{j+1}+c_{1n}}-1
=\frac{\hat{s}^2_{j}-\hat{s}^2_{j+1}}{\hat{s}^2_{j+1}+c_{1n}}=O_p(\max\{C^4_n/c_{1n},1/(nc_{1n})\}).
\end{eqnarray*}
Consider two situations. When $C_n=n^{-\alpha}$ with $1/4\leq\alpha<1/2$, one has $O_p(\max\{C^4_n/c_{1n},1/(nc_{1n})\})=O_p(1/(nc_{1n}))$.
As $c_{1n}c_{2n}n \rightarrow \infty$, we derive
\begin{eqnarray*}
\frac{\hat{s}^*_{(j+1)}+c_{2n}}{\hat{s}^*_{j}+c_{2n}}
=\frac{c_{2n}+O_p\{1/(nc_{1n})\}}{c_{2n}+O_p\{1/(nc_{1n})\}}\rightarrow  1 >\tau.
\end{eqnarray*}
When $C_n=n^{-\alpha}$ with $0\leq\alpha<1/4$, one has $O_p(\max\{C^4_n/c_{1n},1/(nc_{1n})\})=O_p(C^4_n/c_{1n})$.
As $c_{1n}c_{2n}/C^4_n \rightarrow \infty$, we derive
\begin{eqnarray*}
\frac{\hat{s}^*_{(j+1)}+c_{2n}}{\hat{s}^*_{j}+c_{2n}}
=\frac{c_{2n}+O_p(C^4_n/c_{1n})\}}{c_{2n}+O_p(C^4_n/c_{1n})\}}\rightarrow  1 >\tau.
\end{eqnarray*}
 Note that the ratio at $j=q_1$ is also a local minimum  converging to zero in probability. But in probability the thresholding criterion prevents this value to be an estimate of $q$. Thus we conclude that as $n\rightarrow \infty$, $P(\hat{q}=q) \rightarrow 1$.
The results of (\ref{star1}) have been proved.
\hfill$\Box$

\textbf{Proof of Theorem~\ref{theorem-factor}. }
As the arguments used is very similar to that for Theorem~\ref{theorem0}, and the results are also very similar to those in (\ref{lambdastar}), (\ref{lambdastar1}) and (\ref{ratio1}), we then only give an outline. Following the similar justifications as those of Theorem~1 in Wang (2012), we have that with a probability tending to one, $\hat{s}_d>\kappa$ with some positive constant $\kappa$ and for any $i>d$, $\hat{\lambda}_i=O_p(1/m)$. More details can be referred to Wang (2012).
Due to the monotonicity of the function $x/(1+x)$
 there is positive constant $C$ such that  for $i\le d$,
\begin{eqnarray*}
C\ge \lim_{n\to \infty}\hat{s}_{i} &= &\lim_{n\to \infty} \frac{\hat{\lambda}_{i}}{\hat{\lambda}_{i}+1}>\frac{\kappa}{\kappa+1}=:c>0,
\end{eqnarray*}
and for any $i>d$,
\begin{eqnarray*}
\hat{s}_{i} &= &\frac{\hat{\lambda}_{i}}{\hat{\lambda}_{i}+1}=O_p(1/m).
\end{eqnarray*}
Since $c_{n,p} \rightarrow 0$, we have
\begin{eqnarray*}
\hat{s}^*_d=\frac{\hat{s}_{d}+c_{n,p}}{\hat{s}_{d+1}+c_{n,p}}-1
>\frac{\kappa/(\kappa+1)+c_{n,p}}{O_p(1/m)+c_{n,p}}-1\rightarrow \infty,
\end{eqnarray*}
and for any $i>d$,
\begin{eqnarray*}
\hat{s}^*_i=\frac{\hat{s}_{i}+c_{n,p}}{\hat{s}_{i+1}+c_{n,p}}-1
=\frac{\hat{s}_{i}-\hat{s}_{i+1}}{\hat{s}_{i}+c_{n,p}}=O_p\{1/(c_{n,p}m)\}.
\end{eqnarray*}
The results can be applied to prove the following:
\begin{eqnarray*}
\lim_{n\rightarrow \infty}\frac{\hat{s}^*_{(i+1)}+\tilde{c}_{n,p}}{\hat{s}^*_{i}+\tilde{c}_{n,p}}=\left\{\begin{array}{ll}
\ge c, & {\rm{\ for}}\   1=i < d-1,\\
+\infty, & {\rm{\ for}}\   i = d-1,\\
0 < \tau, & {\rm{\ for}}\   i = d,\\
1>\tau, & {\rm{\ for}}\  d+1 \leq j \leq p-2.
\end{array}\right.
\end{eqnarray*}
The details are omitted. Therefore, we  conclude that as $n\rightarrow \infty$, $P(\hat{d}=d)=1$. \hfill$\Box$

\

\leftline{\large\bf References}

\begin{description}
\item Ahn, S. C., and Horenstein, A. R. (2013). Eigenvalue ratio test for the number of factors. {\it Econometrica}, {\bf 81}, 1203-1227.

\item Bai, J., and  Ng, S. (2002). Determining the number of factors in approximate factor models. {\it Econometrica}, {\bf 70}, 191-221.

\item Bi, X. and Qu, A. (2015). Sufficient dimension reduction for longitudinal data. {\it Statistica Sinica}, doi:http://dx.doi.org/10.5705/ss.2013.168

\item Bura, E. and Cook, R. D. (2001). Extending sliced inverse regression: The weighted chi-squared test. {\it Journal of the American Statistical Association}, {\bf 96}, 996-1003.

\item Bura, E. and Yang, J. (2011). Dimension estimation in sufficient dimension reduction: a unifying approach. {\it Journal of Multivariate Analysis}, {\bf 102}, 130-142.

\item Caner, M. and  Han, X. (2014). Selecting the correct number of factors in approximate factor models: the large panel case with group bridge estimators. {\it Journal of Business and Economic Statistics}, {\bf 32}, 359-374.

\item Chen, X., Zou, C. and Cook, R. D. (2010). Coordinate-independent sparse sufficient dimension reduction and variable selection. {\it The Annals of Statistics}, {\bf 38}, 3696-3723.

\item Cook, R. D. (1998). {\it Regression graphics: ideas for studying regressions through graphics. } {New York: Wiley.}

\item Cook, R. D., and Li, B. (2002). Dimension reduction for conditional mean in regression. {\it Annals of Statistics}, {\bf 30}, 455-474.

\item Cook, R. D. and Li, B. (2004). Determining the dimension of iterative Hessian transformation. {\it The Annals of Statistics}, {\bf 32}, 2501-2531.
\item Dong, Y. X. and Li, B. (2010). Dimension reduction for non-elliptically distributed predictors: second-order methods. {\it Biometrika}, {\bf 97}, 279-294.
\item Duan, N. and Li, K. C. (1991). Slicing regression: a link-free regression method.  {\it The Annals of Statistics}, {\bf 19}, 505-530.

\item Fan, J., Fan, Y. and  Lv, J. (2008). High dimensional covariance matrix estimation using a factor model. {\it Journal of Econometrics}, {\bf 147}, 186-197.

\item Fan, J., Liu, H. and Wang, W. (2015). Large covariance estimation through elliptical factor models. {\it arXiv preprint arXiv:1507.08377.} \url{http://arxiv.org/pdf/1507.08377.pdf}
\item Guan, Y., Xie, C.L. and Zhu, Lixing (2016). Sufficient dimension reduction with mixture multivariate skew elliptical distributions, {\it Statistica Sinica}, in press

\item Guo, X., Wang, T. and Zhu, L. X. (2015). Model checking for parametric single-index models: a dimension reduction model-adaptive approach. {\it Journal of the Royal Statistical Society: Series B (Statistical Methodology),} online.


\item Jiang, C. R., Yu, W. and Wang, J. L. (2104). Inverse regression for longitudinal data. {\it The Annals of Statistics}, {\bf 42}, 563 - 591.

\item Li, B. and Dong, Y. X. (2009). Dimension reduction for non-elliptically distributed predictors.
{\it The Annals of Statistics} {\bf 37}, 1272-1298.
\item Li, L. and Lu, W. B. (2008). Sufficient dimension reduction with missing predictors. {\it Journal of the American Statistical Association}, {\bf 103}, 822-831.


\item Li, B., Wen, S. Q. and Zhu, L. X. (2008). On a projective resampling method for dimension reduction with multivariate responses. {\it Journal of the American Statistical Association}, {\bf 103}, 1177-1186.

\item Li, K. C. (1991). Sliced inverse regression for dimension reduction. {\it Journal of the American Statistical Association}, {\bf 86}, 316-327.

\item Li, L., Cook, R. D. and Nachtsheim, C. J. (2005). Model-free variable selection. {\it Journal of the Royal Statistical Society: Series B (Statistical Methodology)}, {\bf 67}, 285-299.

\item Li, Y. and Hsing, T. (2010). Deciding the dimension of effective dimension reduction space for functional and high-dimensional data, {\it Annals of Statistiscs},
{\bf 38}, 3028 - 3062.

\item Luo, R., Wang, H. and Tsai, C. L. (2009). Contour projected dimension reduction. {\it The Annals of Statistics}, {\bf 37}, 3743-3778.


\item Ma, Y. and Zhang, X. (2015). A validated information criterion to determine the structural dimension in dimension reduction models. {\it Biometrika}, {\bf 102}, 409-420


\item Schott, J. R. (1994). Determining the dimensionality in sliced inverse regression. {\it Journal of the American Statistical Association}, {\bf 89}, 141-148.

\item Stock, J. and  Watson, M. (2005). Implications of dynamic factor models for var analysis. {\it NBER Working Paper.}

\item Velilla, S. (1998). Assessing the number of linear components in a general regression problem. {\it Journal of the American Statistical Association}, {\bf 93}, 1088-1098.

\item Xia, Q., Xu, W. and Zhu, L. (2015). Consistently determining the number of factors in multivariate volatility modelling. {\it Statistica Sinica}, {\bf 25}, 1025-1044.

\item Xia, Y. (2007). A constructive approach to the estimation of dimension reduction directions. {\it The Annals of Statistics}, {\bf 35}, 2654-2690.

\item Xia, Y. C., Tong, H., Li, W. K. and Zhu, L. X. (2002). An adaptive estimation of dimension reduction space. {\it Journal of the Royal Statistical Society: Series B}, {\bf 64}, 363-410.

\item  Wang, H. (2012). Factor profiled sure independence screening. {\it Biometrika}, {\bf 99}, 15-28.

\item Wang, H. and Xia, Y. (2008).  Sliced regression for dimension reduction. {\it Journal of the American Statistical Association}, {\bf 103}, 811-821.

\item  Wang, T., Xu, P. and Zhu,  L. X. (2015). Variable selection and estimation for semi-parametric multipleindex models. {\it Bernoulli}, {\bf 21}, 242¨C275.

\item Wang, T. and Zhu, L (2013). Sparse sufficient dimension reduction using optimal scoring. {\it Computational Statistics and Data Analysis,} {\bf  57}, 223¨C232.

\item Wang, Q and Yin, X. R. (2008). A nonlinear multi-dimensional variable selection method for high dimensional data: Sparse MAVE. {\it Computational Statistics and Data Analysis}, {\bf 52}, 4512-4520.

\item Wu, Y. and Li, L. (2011). Asymptotic properties of sufficient dimension reduction with a diverging number of predictors. {\it Statistica Sinica}, {\bf 21}, 707-730.

\item Zeng, P. (2008). Determining the dimension of the central subspace and central mean subspace. {\it Biometrika}, {\bf 95}, 469-479.

\item Zhu, L. P., Zhu, L. X. , Ferr\'{e}, L. and Wang, T. (2010). Sufficient dimension reduction through discretization-expectation estimation. {\it Biometrika}, {\bf 97}, 295-304.

\item Zhu, L. P. and Zhu, L. X. (2009). Dimension reduction for conditional variance in regressions. {\it Statistica Sinica}, {\bf 2}, 869-883.

\item Zhu, L. X., Miao, B. Q. and Peng, H. (2006). On sliced inverse regression with high dimensional covariates. {\it Journal of the American Statistical Association}, {\bf 101}, 630-643.

\item Zhu, L. X. and Fang, K. T. (1996). Asymptotics for the kernel estimates of sliced inverse regression. {\it The Annals of Statistics}, {\bf 24}, 1053-1067.

\item Zhu, L. X. and Ng, K. W. (1995). Asymptotics for sliced inverse regression. {\it Statistica Sinica}, {\bf 5}, 727-736.


\item Zhu, X., Guo, X. and Zhu, L. (2016). An adaptive-to-model test for partially parametric single-index models. {\it Statistics and Computing.} In press. \url{http://link.springer.com/article/10.1007/s11222-016-9680-z}
\end{description}

\newpage

\begin{table}[htb!]\caption{The frequencies of estimated dimension  for Example~1, $q=3$.   \label{model1-1}
\vspace{-0.05cm}}
\centering
 {\tiny\scriptsize\hspace{12.5cm}
\renewcommand{\arraystretch}{1}\tabcolsep 0.1cm
\begin{tabular}{l|ccccc|ccccc}
\hline
$n$&\multicolumn{5}{c}{$200$}&\multicolumn{5}{c}{$800$}
\\ \hline
$p=5$ &$\hat{q}=1$ &$\hat{q}=2$ &$\hat{q}=3$& $\hat{q}=4$ &$\hat{q}=5$ &$\hat{q}=1$ &$\hat{q}=2$ &$\hat{q}=3$& $\hat{q}=4$ &$\hat{q}=5$\\
$\rm{TDRR}$&0.210 &0.628 &0.162 &0 &0 &0.014    &0.280  &0.706 &0 &0 \\
RRE        &0.986 &0.014 &0.030 &0 &0 &1.000    &0      &0     &0 &0 \\
RE         &0.844 &0.088 &0.068 &0 &0 &0.710    &0.012  &0.278 &0 &0 \\
BIC        &1.000 &0     &0     &0 &0 &1.000    &0      &0     &0 &0 \\
\hline
$p=10$ &$\hat{q}=1$ &$\hat{q}=2$ &$\hat{q}=3$& $\hat{q}=4$ &$\hat{q}=5$ &$\hat{q}=1$ &$\hat{q}=2$ &$\hat{q}=3$& $\hat{q}=4$ &$\hat{q}=5$\\
$\rm{TDRR}$&0.164  &0.542   &0.292   &0.002  &0  &0.006    &0.254    &0.740    &0    &0  \\
RRE        &0.992  &0.008   &0.004   &0      &0  &1.000    &0        &0        &0    &0  \\
RE         &0.932  &0.062   &0.006   &0      &0  &0.980    &0.004    &0.016    &0    &0  \\
BIC        &0.962  &0.038   &0       &0      &0  &0.994    &0.006    &0        &0    &0  \\
\hline
$p=20$ &$\hat{q}=1$ &$\hat{q}=2$ &$\hat{q}=3$& $\hat{q}=4$ &$\hat{q}=5$ &$\hat{q}=1$ &$\hat{q}=2$ &$\hat{q}=3$& $\hat{q}=4$ &$\hat{q}=5$\\
$\rm{TDRR}$&0.104   &0.464   &0.408    &0.024   &0  &0.010     &0.240    &0.746   &0.004    &0  \\
RRE        &0.990   &0.010   &0.004    &0       &0  &0.998     &0.002    &0       &0        &0  \\
RE         &0.970   &0.026   &0.004    &0       &0  &0.996     &0.004    &0       &0        &0  \\
BIC        &0.414   &0.586   &0        &0       &0  &0.442     &0.558    &0       &0        &0  \\
\hline
$p=30$ &$\hat{q}=1$ &$\hat{q}=2$ &$\hat{q}=3$& $\hat{q}=4$ &$\hat{q}=5$ &$\hat{q}=1$ &$\hat{q}=2$ &$\hat{q}=3$& $\hat{q}=4$ &$\hat{q}=5$\\                                     $\rm{TDRR}$&0.072    &0.330    &0.436   &0.162  &0  & 0.026&0.244   &0.690   &0.040 &0    \\
RRE        &0.990    &0.010    &0       &0      &0  &1.000    &0       &0       &0     &0    \\
RE         &0.978    &0.022    &0       &0      &0  &1.000    &0       &0       &0     &0    \\
BIC        &0.026    &0.952    &0.022   &0      &0  &0.040    &0.948   &0.012   &0     &0    \\
\hline
$p=40$ &$\hat{q}=1$ &$\hat{q}=2$ &$\hat{q}=3$& $\hat{q}=4$ &$\hat{q}=5$ &$\hat{q}=1$ &$\hat{q}=2$ &$\hat{q}=3$& $\hat{q}=4$ &$\hat{q}=5$\\
$\rm{TDRR}$&0.044  &0.250 &0.418 &0.278 &0.010 &0.026      &0.246    &0.592   &0.136 &0  \\
RRE        &0.990  &0.010 &0     &0     &0     &1.000      &0        &0       &0     &0  \\
RE         &0.980  &0.020 &0     &0     &0     &1.000      &0        &0       &0     &0  \\
BIC        &0      &0.710 &0.290 &0     &0     &0.002      &0.874    &0.124   &0     &0  \\
\hline
\end{tabular}
}\end{table}

\begin{table}[htb!]\caption{Proportion of correct decisions about the structure dimension $q$ when SIR is used in  Example 2.  \label{Example2}
\vspace{0.25cm}}
\centering
 {\tiny\scriptsize\hspace{12.5cm}
\renewcommand{\arraystretch}{1}\tabcolsep 0.2cm
\begin{tabular}{cccc|ccc}
\hline
$n$&\multicolumn{3}{c}{$400$}&\multicolumn{3}{c}{$800$}
\\
\hline
$p$ &BIC &ST & $\rm{TDRR}$&BIC  &ST &$\rm{TDRR}$\\
\hline
10&0.78&0.39&0.852&0.94&0.70&0.986\\
20&0.73&0.24&0.692&0.85&0.46&0.906\\
30&0.69&0.17&0.560&0.81&0.31&0.826 \\
40&0.43&0.11&0.250&0.63&0.27&0.738\\
\hline
\end{tabular}
}
\end{table}

%
%
%
%
%
%
%
%

\begin{table}[htb!]\caption{The frequencies of estimated dimensions for Examples~3 and 4, $a=1/n^{1/4}$, $q_1=1$, $q=2$.   \label{Example4}
\vspace{0.05cm}}
\centering
 {\tiny\scriptsize\hspace{12.5cm}
\renewcommand{\arraystretch}{1}\tabcolsep 0.1cm
\begin{tabular}{cc|ccc|ccc|ccc|ccc}
\hline
&&\multicolumn{3}{c}{$\rm{TDRR}$}
& \multicolumn{3}{c}{RRE}
& \multicolumn{3}{c}{RE }
& \multicolumn{3}{c}{BIC}\\
\hline
Example 3&p&$\hat{q}=1$ &$\hat{q}=2$ &$\hat{q}=3$
&$\hat{q}=1$ &$\hat{q}=2$ &$\hat{q}=3$
&$\hat{q}=1$ &$\hat{q}=2$ &$\hat{q}=3$
&$\hat{q}=1$ &$\hat{q}=2$ &$\hat{q}=3$\\
\hline
&5    &0.220&0.780&0    &0.988&0.012&0&0.812&0.186&0.002&1.000&0    &0\\
$n=200$
&7    &0.214&0.786&0    &0.994&0.006&0&0.842&0.154&0.002&1.000&0    &0\\
&10   &0.180&0.802&0.018&0.980&0.020&0&0.910&0.090&0    &0.962&0.038&0\\
\hline
&5    &0.086&0.914&0    &1.000&0    &0&0.816&0.184&0    &1.000&0    &0\\
$n=800$
&7    &0.050&0.950&0    &1.000&0    &0&0.890&0.110&0    &1.000&0    &0\\
&10   &0.070&0.930&0    &1.000&0    &0&0.960&0.040&0    &0.994&0.006&0\\
\hline
Example 4&p&$\hat{q}=1$ &$\hat{q}=2$ &$\hat{q}=3$
&$\hat{q}=1$ &$\hat{q}=2$ &$\hat{q}=3$
&$\hat{q}=1$ &$\hat{q}=2$ &$\hat{q}=3$
&$\hat{q}=1$ &$\hat{q}=2$ &$\hat{q}=3$\\
\hline
&5    &0.082&0.918&0    &0.786&0.214&0    &0.386&0.594&0.020&0.964&0.036&0\\
$n=200$
&7    &0.058&0.936&0.006&0.812&0.188&0    &0.498&0.482&0.020&0.872&0.128&0\\
&10   &0.042&0.942&0.016&0.780&0.220&0    &0.552&0.442&0.006&0.590&0.410&0\\
\hline
&5    &0.002&0.998&0    &0.744&0.256&0    &0.124&0.870&0.006&0.934&0.066&0\\
$n=800$
&7    &0    &1.000&0    &0.814&0.186&0    &0.164&0.836&0    &0.796&0.204&0\\
&10   &0    &1.000&0    &0.840&0.160&0    &0.354&0.646&0    &0.462&0.538&0\\
\hline
\end{tabular}
}
\end{table}

\begin{table}[htb!]\caption{\small The frequencies of estimation for the number of common factors with 
$f_t\sim N(0, \Sigma_i)$ in Example~5. \label{factor-N}
\vspace{-0.5cm}}
\centering
 {\tiny\scriptsize\hspace{12.5cm}
\renewcommand{\arraystretch}{1}\tabcolsep 0.15cm
\begin{tabular}{cc|ccc|ccc|ccc|ccc}
\hline
&&\multicolumn{3}{c}{TDRR}
& \multicolumn{3}{c}{RRE}
& \multicolumn{3}{c}{RE}
& \multicolumn{3}{c}{BIC}\\
\hline
$f_t$&$p/n$ &$\hat{d}<4$ &$\hat{d}=4$&$\hat{d}>4$ &$\hat{d}<4$ &$\hat{d}=4$&$\hat{d}>4$
&$\hat{d}<4$ &$\hat{d}=4$&$\hat{d}>4$ &$\hat{d}<4$ &$\hat{d}=4$&$\hat{d}>4$\\
\hline
$N(0,\Sigma_1)$
&1    &0   &0.996&0.004&0    &1.000&0    &0    &1.000&0    &0.004&0.996    &0\\
n=50
&2    &0   &1.000&0    &0    &1.000&0    &0    &1.000&0    &0    &1.000    &0\\
&4    &0   &0.936&0.064&0    &1.000&0    &0    &1.000&0    &0    &1.000    &0\\
 \hline
$N(0,\Sigma_1)$
&1    &0   &1.000&0    &0    &1.000&0    &0    &1.000&0    &0    &1.000    &0\\
n=100
&2    &0   &1.000&0    &0    &1.000&0    &0    &1.000&0    &0    &1.000    &0\\
&4    &0   &1.000&0    &0    &1.000&0    &0    &1.000&0    &0    &1.000    &0\\
 \hline
$N(0,\Sigma_1)$
&1    &0   &1.000&0    &0    &1.000&0    &0    &1.000&0    &0    &1.000    &0\\
n=200
&2    &0   &1.000&0    &0    &1.000&0    &0    &1.000&0    &0    &1.000    &0\\
&4    &0   &1.000&0    &0    &1.000&0    &0    &1.000&0    &0    &1.000    &0\\
\hline
$N(0,\Sigma_2)$
&1    &0.006&0.980&0.014&0.944&0.056&0    &0.908&0.092&0    &0.992&0.008&0\\
n=50
&2    &0.002&0.996&0.002&0.988&0.012&0    &0.978&0.022&0    &0.942&0.058&0\\
&4    &0    &0.926&0.074&0.996&0.004&0    &0.988&0.012&0    &0.940&0.060&0\\
 \hline
$N(0,\Sigma_2)$
&1    &0    &1.000&0    &0.998&0.002&0    &0.992&0.008&0    &0.998&0.002&0\\
n=100
&2    &0    &1.000&0    &0.874&0.126&0    &0.728&0.272&0    &0.376&0.624&0\\
&4    &0    &1.000&0    &0.642&0.358&0    &0.362&0.638&0    &0.002&0.998&0\\
 \hline
$N(0,\Sigma_2)$
&1    &0    &1.000&0    &0.422&0.578&0    &0.204&0.796&0    &0.020&0.980&0\\
n=200
&2    &0    &1.000&0    &0.016&0.984&0    &0.002&0.998&0    &0    &1.000&0\\
&4    &0    &1.000&0    &0.002&0.998&0    &0    &0.998&0    &0    &1.000&0\\
\hline
$N(0,\Sigma_3)$
&1    &0.002&0.992&0.006&1.000&0    &0    &1.000&0    &0    &1.000&0    &0\\
n=50
&2    &0    &1.000&0    &0.992&0.008&0    &0.976&0.024&0    &0.842&0.158&0\\
&4    &0    &0.924&0.076&0.986&0.014&0    &0.970&0.030&0    &0.170&0.830&0\\
 \hline
$N(0,\Sigma_3)$
&1    &0    &1.000&0   &0.918&0.082&0    &0.836&0.164&0    &0.060&0.940&0\\
n=100
&2    &0    &1.000&0   &0.412&0.588&0    &0.212&0.788&0    &0    &1.000&0\\
&4    &0    &1.000&0   &0.378&0.622&0    &0.146&0.854&0    &0    &1.000&0\\
 \hline
$N(0,\Sigma_3)$
&1    &0    &1.000&0   &0.018&0.982&0    &0.002&0.998&0    &0    &1.000&0\\
n=200
&2    &0    &1.000&0   &0    &1.000&0    &0    &1.000&0    &0    &1.000&0\\
&4    &0    &1.000&0   &0    &1.000&0    &0    &1.000&0    &0    &1.000&0\\
\hline
\end{tabular}
}
\end{table}

\begin{table}[htb!]\caption{\small The frequencies of estimation for the number of common factors
with $f_t\sim t(2.5, \Sigma_i)$ 
in Example~5. \label{factor-T}
\vspace{-0.5cm}}
\centering
 {\tiny\scriptsize\hspace{12.5cm}
\renewcommand{\arraystretch}{1}\tabcolsep 0.15cm
\begin{tabular}{cc|ccc|ccc|ccc|ccc}
\hline
&&\multicolumn{3}{c}{TDRR}
& \multicolumn{3}{c}{RRE}
& \multicolumn{3}{c}{RE}
& \multicolumn{3}{c}{BIC}\\
\hline
$f_t$&$p/n$ &$\hat{d}<4$ &$\hat{d}=4$&$\hat{d}>4$ &$\hat{d}<4$ &$\hat{d}=4$&$\hat{d}>4$
&$\hat{d}<4$ &$\hat{d}=4$&$\hat{d}>4$ &$\hat{d}<4$ &$\hat{d}=4$&$\hat{d}>4$\\
\hline
$t(2.5,\Sigma_1)$
&1    &0   &0.992&0.008&0.012&0.988&0    &0.010&0.990&0    &0.096&0.904&0\\
n=50
&2    &0   &1.000&0    &0.010&0.990&0    &0.010&0.990&0    &0.034&0.966&0\\
&4    &0   &0.910&0.090&0.016&0.984&0    &0.010&0.990&0    &0.016&0.984&0\\
 \hline
$t(2.5,\Sigma_1)$
&1    &0   &1.000&0    &0.006&0.994&0    &0.004&0.996&0    &0.018&0.982    &0\\
n=100
&2    &0   &1.000&0    &0    &1.000&0    &0    &1.000&0    &0    &1.000    &0\\
&4    &0   &1.000&0    &0.002&0.998&0    &0.002&0.998&0    &0    &1.000    &0\\
 \hline
$t(2.5,\Sigma_1)$
&1    &0   &1.000&0    &0    &1.000&0    &0    &1.000&0    &0    &1.000    &0\\
n=200
&2    &0   &1.000&0    &0    &1.000&0    &0    &1.000&0    &0    &1.000    &0\\
&4    &0   &1.000&0    &0    &1.000&0    &0    &1.000&0    &0    &1.000    &0\\
\hline
$t(2.5,\Sigma_2)$
&1    &0   &0.984&0.016&0.794&0.206&0    &0.742&0.258&0    &0.994&0.006    &0\\
n=50
&2    &0   &1.000&0    &0.568&0.432&0    &0.504&0.496&0    &0.900&0.100    &0\\
&4    &0   &0.906&0.094&0.628&0.372&0    &0.508&0.492&0    &0.846&0.154    &0\\
 \hline
$t(2.5,\Sigma_2)$
&1    &0    &1.000&0    &0.564&0.436&0    &0.510&0.490&0    &0.910&0.090    &0\\
n=100
&2    &0    &1.000&0    &0.218&0.782&0    &0.156&0.844&0    &0.350&0.650    &0\\
&4    &0    &1.000&0    &0.094&0.906&0    &0.066&0.934&0    &0.080&0.920    &0\\
 \hline
$t(2.5,\Sigma_2)$
&1    &0    &1.000&0    &0.032&0.968&0    &0.020&0.980&0    &0.068&0.932&0\\
n=200
&2    &0    &1.000&0    &0.012&0.988&0    &0.012&0.988&0    &0.010&0.990&0\\
&4    &0    &1.000&0    &0.020&0.980&0    &0.018&0.982&0    &0.012&0.988&0\\
\hline
$t(2.5,\Sigma_3)$
&1    &0    &0.988&0.012&0.708&0.292&0    &0.666&0.334&0    &0.972&0.028&0\\
n=50
&2    &0    &1.000&0    &0.446&0.554&0    &0.384&0.616&0    &0.642&0.358&0\\
&4    &0    &0.906&0.094&0.186&0.814&0    &0.146&0.854&0    &0.144&0.856&0\\
 \hline
$t(2.5,\Sigma_3)$
&1    &0    &1.000&0   &0.234&0.766&0    &0.190&0.810&0    &0.276&0.724&0\\
n=100
&2    &0    &1.000&0   &0.100&0.900&0    &0.076&0.924&0    &0.056&0.944&0\\
&4    &0    &1.000&0   &0.024&0.976&0    &0.016&0.984&0    &0.014&0.986&0\\
 \hline
$t(2.5,\Sigma_3)$
&1    &0    &1.000&0   &0.014&0.986&0    &0.014&0.986&0    &0.014&0.986&0\\
n=200
&2    &0    &1.000&0   &0.012&0.988&0    &0.010&0.990&0    &0.002&0.998&0\\
&4    &0    &1.000&0   &0.006&0.994&0    &0.002&0.998&0    &0    &1.000&0\\
\hline
\end{tabular}
}
\end{table}

\newpage

\begin{figure}
  \centering
  \includegraphics[width=7.250cm]{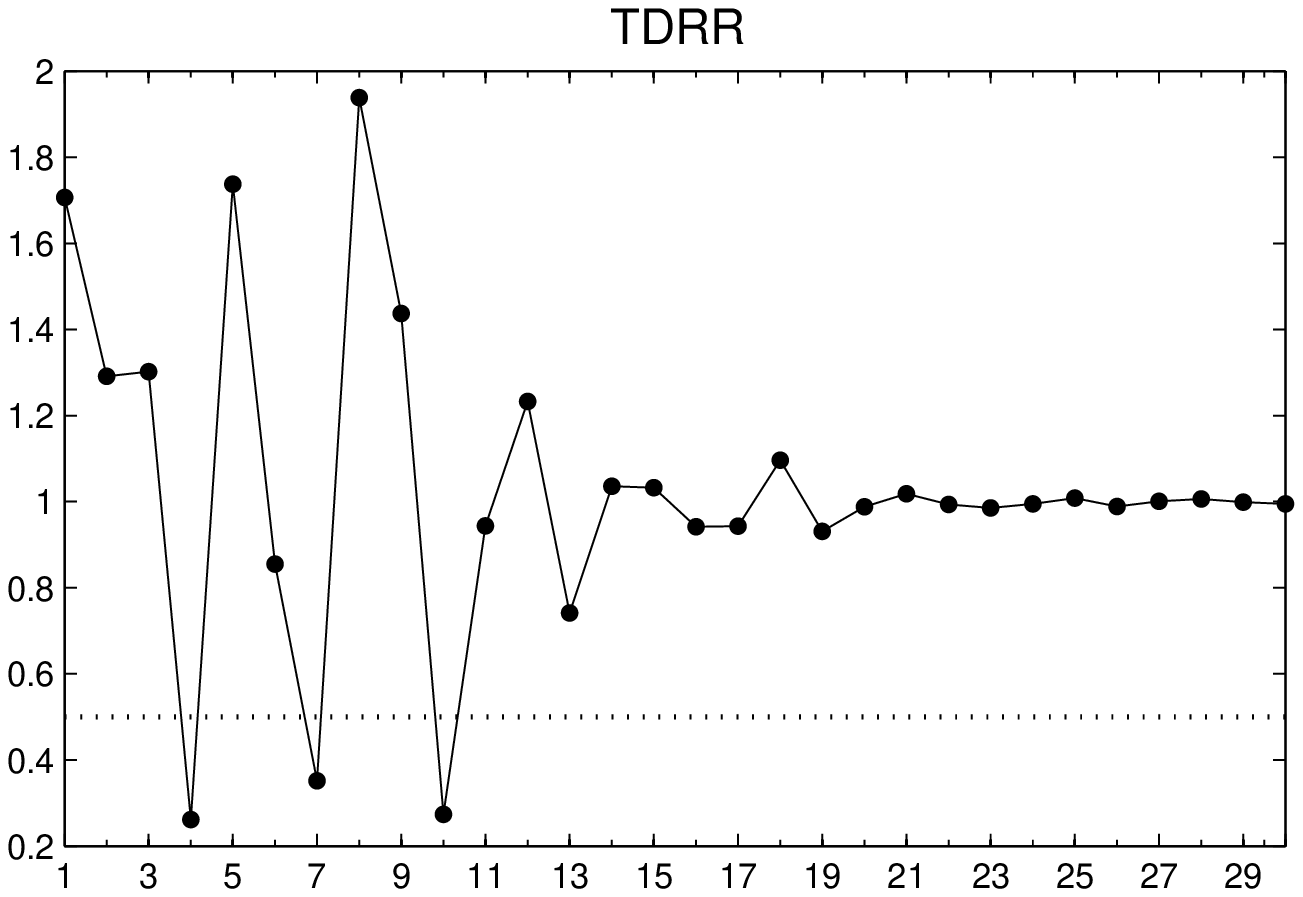}
  \includegraphics[width=7.250cm]{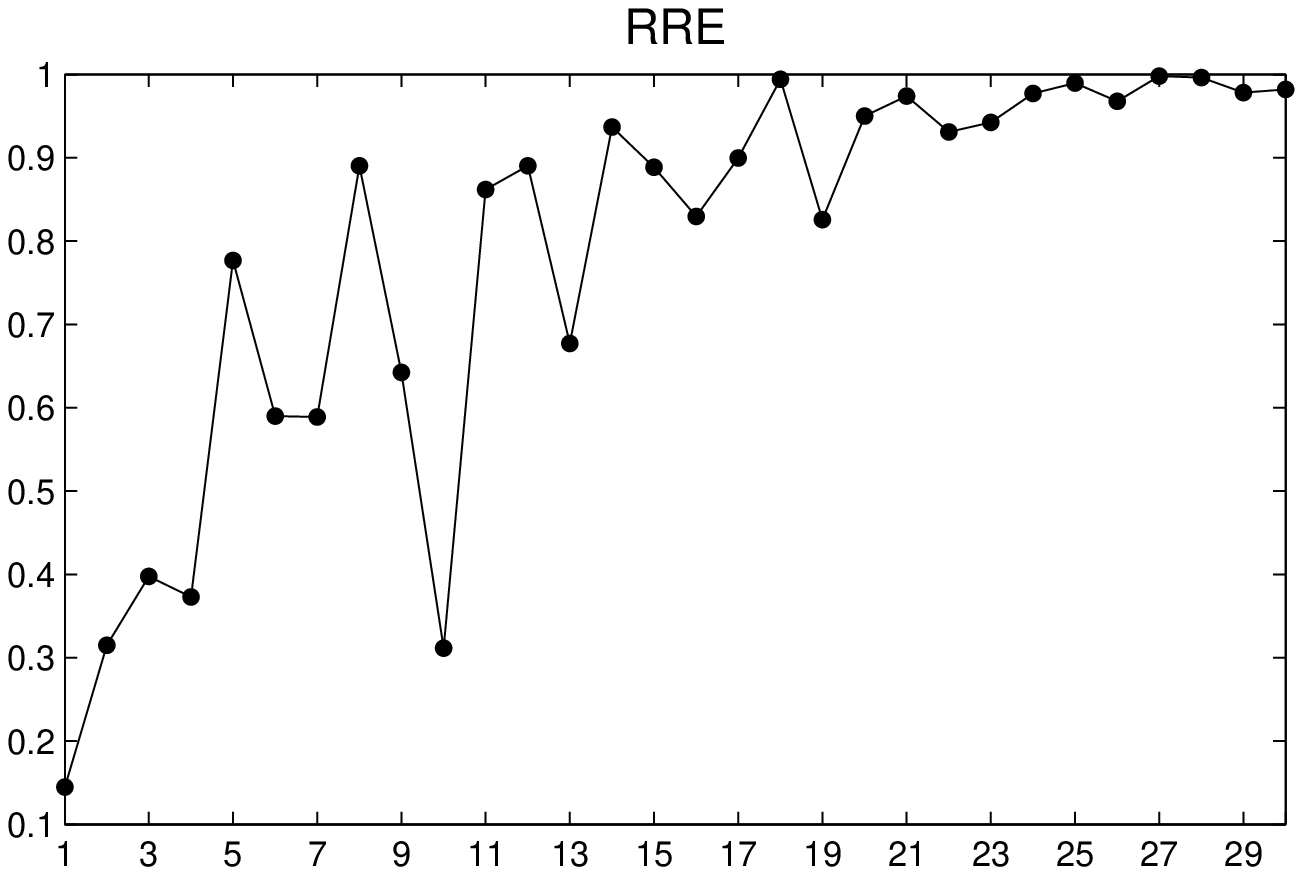}
  \includegraphics[width=7.250cm]{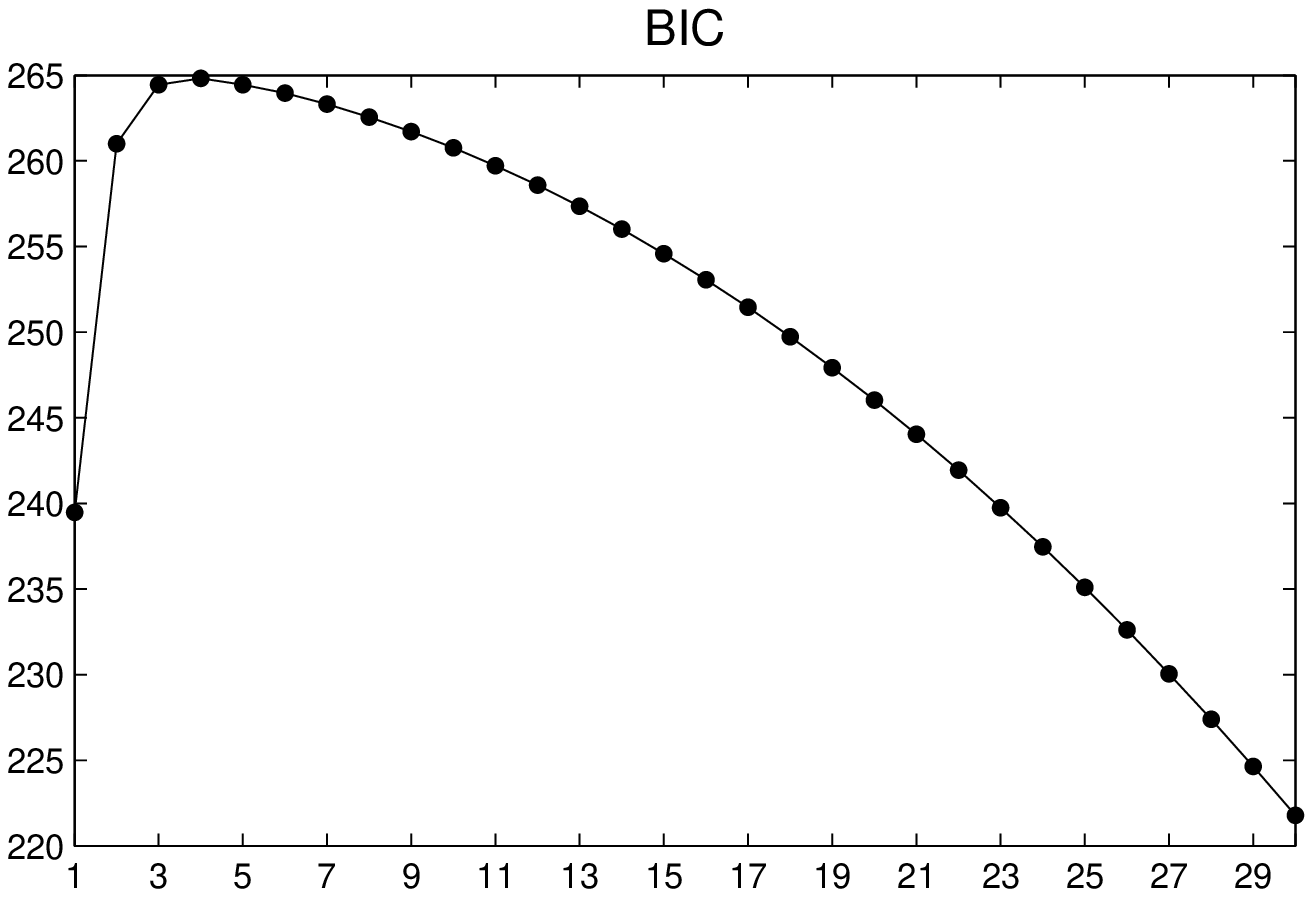}
  \includegraphics[width=7.250cm]{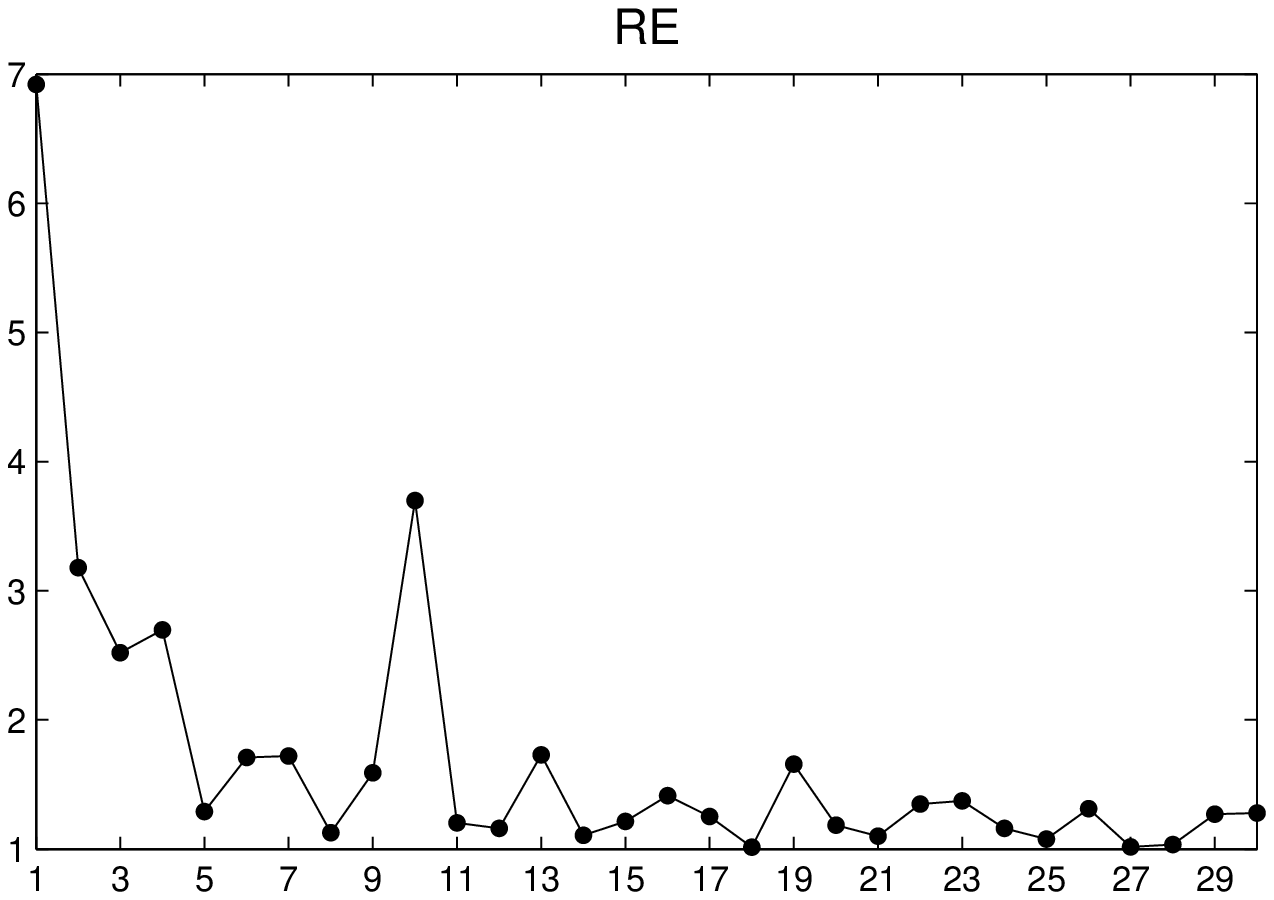}\\
  \caption{Plots of the first $30$ components for $\{\hat{s}^*_{j+1}+\frac{\log{(m)}}{5\sqrt{m}}\}/\{\hat{s}^*_j+\frac{\log{(m)}}{5\sqrt{m}}\}$ (TDRR),  $\{\hat{\lambda}_j+\frac{\log{n}}{10n}\}/\{\hat{\lambda}_{j+1}+\frac{\log{n}}{10n}\}$ (RRE),
 \{$\frac{\hat{\lambda}_j}{\hat{\lambda}_{j+1}}$\} (RE) and
  $\{\frac{n\sum_{l=1}^j\{\log(1+\hat{\lambda}_l)+\hat{\lambda}_l\}}
  {2\sum_{l=1}^p\{\log(1+\hat{\lambda}_l)+\hat{\lambda}_l\}}-\log(n)\frac{j(j+1)}{p}\}$ (BIC), for $j=1,\ldots,30$,  about the approximate factor model of the real data.
  }\label{real-factor}
\end{figure}

\end{document}